\documentclass[final,3p]{elsarticle}
 \usepackage{graphics}
 \usepackage{graphicx}
 \usepackage{epsfig}
\usepackage{amssymb}
 \usepackage{amsthm}
 \usepackage{lineno}
 \usepackage{amsmath}
   \numberwithin{equation}{section}
\usepackage{mathrsfs}

\newtheorem{thm}{Theorem}[section]

\newtheorem{lem}[thm]{Lemma}

\biboptions{sort&compress,square}
\allowdisplaybreaks
\begin{document}
\begin{frontmatter}
\author{Tong Wu}
\ead{wut977@nenu.edu.cn}
\author{Yong Wang\corref{cor3}}
\ead{wangy581@nenu.edu.cn}
\cortext[cor3]{Corresponding author.}

\address{School of Mathematics and Statistics, Northeast Normal University,
Changchun, 130024, China}
\title{The generalized noncommutative residue and the Kastler-Kalau-Walze type theorem }
\begin{abstract}
In this paper, we define the generalized noncommutative residue of the Dirac operator.  And we give the proof of Kastler-Kalau-Walze type theorems for the generalized noncommutative residue on 4-dimensional and 6-dimensional compact manifolds with (resp.without) boundary.
\end{abstract}
\end{frontmatter}
\hskip -0.4 true cm \textit{2010 Mathematics Subject Classification:}
53C40; 53C42.
\section{Introduction}
 Until now, many geometers have studied noncommutative residues. In \cite{Gu,Wo}, authors found noncommutative residues are of great importance to the study of noncommutative geometry. In \cite{Co1}, Connes used the noncommutative residue to derive a conformal 4-dimensional Polyakov action analogy. Connes showed us that the noncommutative residue on a compact manifold $M$ coincided with the Dixmier's trace on pseudodifferential operators of order $-{\rm {dim}}M$ in \cite{Co2}.
And Connes claimed the noncommutative residue of the square of the inverse of the Dirac operator was proportioned to the Einstein-Hilbert action.  Kastler \cite{Ka} gave a
brute-force proof of this theorem. Kalau and Walze proved this theorem in the normal coordinates system simultaneously in \cite{KW} .
Ackermann proved that
the Wodzicki residue  of the square of the inverse of the Dirac operator ${\rm  Wres}(D^{-2})$ in turn is essentially the second coefficient
of the heat kernel expansion of $D^{2}$ in \cite{Ac}.

On the other hand, Wang generalized the Connes' results to the case of manifolds with boundary in \cite{Wa1,Wa2},
and proved the Kastler-Kalau-Walze type theorem for the Dirac operator and the signature operator on lower-dimensional manifolds
with boundary \cite{Wa3}. In \cite{Wa3,Wa4}, Wang computed $\widetilde{{\rm Wres}}[\pi^+D^{-1}\circ\pi^+D^{-1}]$ and $\widetilde{{\rm Wres}}[\pi^+D^{-2}\circ\pi^+D^{-2}]$, where the two operators are symmetric, in these cases the boundary term vanished. But for $\widetilde{{\rm Wres}}[\pi^+D^{-1}\circ\pi^+D^{-3}]$, Wang got a nonvanishing boundary term \cite{Wa5}, and give a theoretical explanation for gravitational action on boundary. In others words, Wang provides a kind of method to study the Kastler-Kalau-Walze type theorem for manifolds with boundary. In \cite{Mh}, Moscovici computed ${\rm  tr}[c(\omega)e^{-tD^2}]$ for a form $\omega$ and the Dirac operator $D$, and used this formula to derive information about the asymptotic distribution of multiplicities in the quasi-regular representation of semisimple Lie grroup modulo a compact discrete subgroup. In \cite{Mj}, Mickelsson and Paycha computed the index of the Dirac operator by the super noncommutative residue. Motivated by \cite{Mh} and \cite{Mj}, we define the generalized noncommutative residue for the Dirac operator which is the generalization of the noncommutative residue and the super noncommutative residue.\\
\indent The motivation of this paper is
to compute  the generalized noncommutative residue $\widetilde{{\rm Wres}}[\pi^+(LD^{-1})\circ\pi^+D^{-1}]$ and $\widetilde{{\rm Wres}}[\pi^+(LD^{-2})\circ\pi^+D^{-2}]$, where $L=c(X_1)c(X_2)c(X_3)\cdot\cdot\cdot c(X_l),$ and prove the Kastler-Kalau-Walze type theorem for the generalized noncommutative
residue on 4-dimensional and 6-dimensional compact manifolds. \\
\indent The paper is organized in the following way. In Section \ref{section:2}, we define the generalized noncommutative residue and get the Kastler-Kalau-Walze type theorem for the generalized noncommutative residue on manifolds without boundary. In Section \ref{section:3} and in Section \ref{section:4},
 we prove the Kastler-Kalau-Walze type theorem for the generalized noncommutative residue on 4-dimensional and 6-dimensional manifolds with boundary respectively.
\section{The Dirac operator and its generalized noncommutative residue}
\label{section:2}
Firstly we recall the definition of Dirac operator. Let $M$ be an $n$-dimensional ($n\geq 3$) oriented compact Riemannian manifold with a Riemannian metric $g^{M}$ and let $\nabla^L$ be the Levi-Civita connection about $g^{M}$. In the fixed orthonormal frame $\{e_1,\cdots,e_n\}$, the connection matrix $(\omega_{s,t})$ is defined by
\begin{equation}
\label{a2}
\nabla^L(e_1,\cdots,e_n)= (e_1,\cdots,e_n)(\omega_{s,t}).
\end{equation}
\indent Let $\epsilon (e_j^*)$,~$\iota (e_j^*)$ be the exterior and interior multiplications respectively, where $e_j^*=g^{TM}(e_j,\cdot)$.
Write
\begin{equation}
\label{a3}
\widehat{c}(e_j)=\epsilon (e_j^* )+\iota
(e_j^*);~~
c(e_j)=\epsilon (e_j^* )-\iota (e_j^* ),
\end{equation}
which satisfies
\begin{align}
\label{a4}
&\widehat{c}(e_i)\widehat{c}(e_j)+\widehat{c}(e_j)\widehat{c}(e_i)=2g^{M}(e_i,e_j);~~\nonumber\\
&c(e_i)c(e_j)+c(e_j)c(e_i)=-2g^{M}(e_i,e_j);~~\nonumber\\
&c(e_i)\widehat{c}(e_j)+\widehat{c}(e_j)c(e_i)=0.\nonumber\\
\end{align}
By \cite{Y}, we have the Dirac operator
\begin{align}
\label{a5}
&D=\sum^n_{i=1}c(e_i)[e_i-\frac{1}{4}\sum_{s,t}\omega_{s,t}
(e_i)c(e_s)c(e_t)].\nonumber\\
\end{align}
\indent Set $L=c(X_1)c(X_2)c(X_3)\cdot\cdot\cdot c(X_l)$, where $X_j=\sum_{\alpha=1}^na_{j\alpha}e_\alpha,~1\leq j\leq l,$ is a smooth vector field. We define the generalized noncommutative residue of $(D^2)^{-\frac{n-2}{2}}$ by ${\rm Wres}[L(D^2)^{-\frac{n-2}{2}}]$. When $L=c(e_1)\cdot\cdot\cdot c(e_n),$ we get the super noncommutative residue of $(D^2)^{-\frac{n-2}{2}}.$ When $X_1=X_2,$ and $l=2,~|X_1|=1,$ we get $-{\rm Wres}[(D^2)^{-\frac{n-2}{2}}].$ \\
By (4.30) in \cite{KW}, we have \\
\begin{align}\label{n}
\int_{|\xi|=1}\sigma_{-n}(D^{-n+2})(x,\xi)d\xi=\frac{(n-2)(4\pi)^{\frac{n}{2}}}{(\frac{n}{2}-1)!}(-\frac{1}{12}s),\nonumber\\
\end{align}
where $s$ is the scalar curvature. Then we have the following theorem
\begin{thm}\label{thm2} If $M$ is an $n$-dimensional compact oriented manifolds without boundary, and $n$ is even, then we get the following equality:
\begin{align}
\label{a29}
{\rm Wres}[L(D^2)^{-\frac{n-2}{2}}]
&=\frac{(n-2)(4\pi)^{\frac{n}{2}}}{(\frac{n}{2}-1)!}\int_{M}\bigg(
-\frac{1}{12}{\rm tr}[c(X_1)c(X_2)\cdot\cdot\cdot c(X_l)]s\bigg)d{\rm Vol_{M}}.
\end{align}
\end{thm}
\indent Now, we need to compute ${\rm tr}[c(X_1)c(X_2)\cdot\cdot\cdot c(X_l)]$ by case. Obviously, we know when $l$ is odd, ${\rm tr}[c(X_1)c(X_2)\cdot\cdot\cdot c(X_l)]$=0. In the following, we compute ${\rm tr}[c(X_1)c(X_2)]$, ${\rm tr}[c(X_1)c(X_2)c(X_3)c(X_4)]$ and ${\rm tr}[c(X_1)c(X_2)c(X_3)c(X_4)c(X_5)c(X_6)]$.\\
$\mathbf{case(1)}$\\
 By $X_j=\sum_{\alpha=1}^na_{j\alpha}e_\alpha$ and (\ref{a4}), we have
\begin{align}
\label{a26}
&{\rm tr}[c(X_1)c(X_2)]=\sum^n_{\alpha,\beta=1}a_{1\alpha}a_{2\beta}{\rm tr}[c(e_\alpha)c(e_\beta)]=-\sum^n_{\alpha=1}a_{1\alpha}a_{2\alpha}{\rm tr}[{\rm \texttt{id}}]=-g(X_1,X_2){\rm tr}[{\rm \texttt{id}}].\nonumber\\
\end{align}
$\mathbf{case(2)}$\\
\begin{align}
\label{a286}
&{\rm tr}[c(X_1)c(X_2)c(X_3)c(X_4)]=\sum^n_{\alpha,\beta,\gamma,\mu=1}a_{1\alpha}a_{2\beta}a_{3\gamma}a_{4\mu}{\rm tr}[c(e_\alpha)c(e_\beta)c(e_\gamma)c(e_\mu)].\nonumber\\
\end{align}
$\mathbf{case(2-a)}$ When $\alpha=\beta, \gamma=\mu$.\\
By (\ref{a4}), we have
\begin{align}
\label{d1}
&\sum^n_{\alpha,\beta,\gamma,\mu=1}a_{1\alpha}a_{2\beta}a_{3\gamma}a_{4\mu}{\rm tr}[c(e_\alpha)c(e_\beta)c(e_\gamma)c(e_\mu)]=\sum^n_{\alpha,\gamma=1}a_{1\alpha}a_{2\alpha}a_{3\gamma}a_{4\gamma}{\rm tr}[{\rm \texttt{id}}].\nonumber\\
\end{align}
$\mathbf{case(2-b)}$ When $\alpha\neq \beta, \alpha=\gamma, \beta=\mu$.\\
By (\ref{a4}), we have
\begin{align}
\label{d2}
&\sum^n_{\alpha,\beta,\gamma,\mu=1}a_{1\alpha}a_{2\beta}a_{3\gamma}a_{4\mu}{\rm tr}[c(e_\alpha)c(e_\beta)c(e_\gamma)c(e_\mu)]=-\sum^n_{\alpha\neq\beta=1}a_{1\alpha}a_{2\beta}a_{3\alpha}a_{4\beta}{\rm tr}[{\rm \texttt{id}}].\nonumber\\
\end{align}
$\mathbf{case(2-c)}$ When $\alpha\neq \beta, \alpha=\mu, \beta=\gamma$.\\
By (\ref{a4}), we have
\begin{align}
\label{d7}
&\sum^n_{\alpha,\beta,\gamma,\mu=1}a_{1\alpha}a_{2\beta}a_{3\gamma}a_{4\mu}{\rm tr}[c(e_\alpha)c(e_\beta)c(e_\gamma)c(e_\mu)]=\sum^n_{\alpha\neq\beta=1}a_{1\alpha}a_{2\beta}a_{3\beta}a_{4\alpha}{\rm tr}[{\rm \texttt{id}}].\nonumber\\
\end{align}
$\mathbf{case(2-d)}$ Other cases.\\
By (\ref{a4}), we have
\begin{align}
\label{d8}
&\sum^n_{\alpha,\beta,\gamma,\mu=1}a_{1\alpha}a_{2\beta}a_{3\gamma}a_{4\mu}{\rm tr}[c(e_\alpha)c(e_\beta)c(e_\gamma)c(e_\mu)]=0.\nonumber\\
\end{align}
Therefore
\begin{align}
\label{d45}
{\rm tr}[c(X_1)c(X_2)c(X_3)c(X_4)]
&=(\sum^n_{\alpha,\gamma=1}a_{1\alpha}a_{2\alpha}a_{3\gamma}a_{4\gamma}-\sum^n_{\alpha\neq\beta=1}a_{1\alpha}a_{2\beta}a_{3\alpha}a_{4\beta} +\sum^n_{\alpha\neq\beta=1}a_{1\alpha}a_{2\beta}a_{3\beta}a_{4\alpha}){\rm tr}[{\rm \texttt{id}}]\nonumber\\
&=(\sum^n_{\alpha,\gamma=1}a_{1\alpha}a_{2\alpha}a_{3\gamma}a_{4\gamma}-\sum^n_{\alpha,\beta=1}a_{1\alpha}a_{2\beta}a_{3\alpha}a_{4\beta}+\sum^n_{\alpha=\beta=1}a_{1\alpha}a_{2\beta}a_{3\alpha}a_{4\beta}\nonumber\\
&+\sum^n_{\alpha,\beta=1}a_{1\alpha}a_{2\beta}a_{3\beta}a_{4\alpha}-\sum^n_{\alpha=\beta=1}a_{1\alpha}a_{2\beta}a_{3\alpha}a_{4\beta}){\rm tr}[{\rm \texttt{id}}]\nonumber\\
&=(\sum^n_{\alpha,\gamma=1}a_{1\alpha}a_{2\alpha}a_{3\gamma}a_{4\gamma}-\sum^n_{\alpha,\beta=1}a_{1\alpha}a_{2\beta}a_{3\alpha}a_{4\beta}+\sum^n_{\alpha,\beta=1}a_{1\alpha}a_{2\beta}a_{3\beta}a_{4\alpha}){\rm tr}[{\rm \texttt{id}}]\nonumber\\
&=[g(X_1,X_2)g(X_3,X_4)-g(X_1,X_3)g(X_2,X_4)+g(X_1,X_4)g(X_2,X_3)]{\rm tr}[{\rm \texttt{id}}].\nonumber\\
\end{align}
$\mathbf{case(3)}$\\
\begin{align}
\label{w1}
&{\rm tr}[c(X_1)c(X_2)c(X_3)c(X_4)c(X_5)c(X_6)]=\sum^n_{\alpha,\beta,\gamma,\mu,\delta,\nu=1}a_{1\alpha}a_{2\beta}a_{3\gamma}a_{4\mu}a_{5\delta}a_{6\nu}{\rm tr}[c(e_\alpha)c(e_\beta)c(e_\gamma)c(e_\mu)c(e_\delta)c(e_\nu)].\nonumber\\
\end{align}
$\mathbf{case(3-a)}$ When $\alpha=\beta, \gamma=\mu, \delta=\nu$.\\
By (\ref{a4}), we have
\begin{align}
\label{w2}
&\sum^n_{\alpha,\beta,\gamma,\mu,\delta,\nu=1}a_{1\alpha}a_{2\beta}a_{3\gamma}a_{4\mu}a_{5\delta}a_{6\nu}{\rm tr}[c(e_\alpha)c(e_\beta)c(e_\gamma)c(e_\mu)c(e_\delta)c(e_\nu)]=-\sum^n_{\alpha,\gamma,\delta=1}a_{1\alpha}a_{2\alpha}a_{3\gamma}a_{4\gamma}a_{5\delta}a_{6\delta}{\rm tr}[{\rm \texttt{id}}].\nonumber\\
\end{align}
$\mathbf{case(3-b)}$ When $\alpha=\beta, \gamma\neq\mu, \gamma=\delta, \mu=\nu$.\\
By (\ref{a4}), we have
\begin{align}
\label{w3}
&\sum^n_{\alpha,\beta,\gamma,\mu,\delta,\nu=1}a_{1\alpha}a_{2\beta}a_{3\gamma}a_{4\mu}a_{5\delta}a_{6\nu}{\rm tr}[c(e_\alpha)c(e_\beta)c(e_\gamma)c(e_\mu)c(e_\delta)c(e_\nu)]=\sum^n_{\alpha,\gamma\neq\mu=1}a_{1\alpha}a_{2\alpha}a_{3\gamma}a_{4\mu}a_{5\gamma}a_{6\mu}{\rm tr}[{\rm \texttt{id}}].\nonumber\\
\end{align}
$\mathbf{case(3-c)}$ When $\alpha=\beta, \gamma\neq\mu, \gamma=\nu, \delta=\mu$.\\
By (\ref{a4}), we have
\begin{align}
\label{w4}
&\sum^n_{\alpha,\beta,\gamma,\mu,\delta,\nu=1}a_{1\alpha}a_{2\beta}a_{3\gamma}a_{4\mu}a_{5\delta}a_{6\nu}{\rm tr}[c(e_\alpha)c(e_\beta)c(e_\gamma)c(e_\mu)c(e_\delta)c(e_\nu)]=-\sum^n_{\alpha,\gamma\neq\mu=1}a_{1\alpha}a_{2\alpha}a_{3\gamma}a_{4\mu}a_{5\mu}a_{6\gamma}{\rm tr}[{\rm \texttt{id}}].\nonumber\\
\end{align}
$\mathbf{case(3-d)}$ When $\alpha\neq \beta, \alpha=\gamma, \beta=\mu, \delta=\nu$.\\
By (\ref{a4}), we have
\begin{align}
\label{w5}
&\sum^n_{\alpha,\beta,\gamma,\mu,\delta,\nu=1}a_{1\alpha}a_{2\beta}a_{3\gamma}a_{4\mu}a_{5\delta}a_{6\nu}{\rm tr}[c(e_\alpha)c(e_\beta)c(e_\gamma)c(e_\mu)c(e_\delta)c(e_\nu)]=\sum^n_{\alpha\neq\beta,\delta=1}a_{1\alpha}a_{2\beta}a_{3\alpha}a_{4\beta}a_{5\delta}a_{6\delta}{\rm tr}[{\rm \texttt{id}}].\nonumber\\
\end{align}
$\mathbf{case(3-e)}$ When $\alpha\neq \beta, \alpha=\gamma, \beta\neq\mu, \beta=\delta, \mu=\nu$.\\
By (\ref{a4}), we have
\begin{align}
\label{w6}
&\sum^n_{\alpha,\beta,\gamma,\mu,\delta,\nu=1}a_{1\alpha}a_{2\beta}a_{3\gamma}a_{4\mu}a_{5\delta}a_{6\nu}{\rm tr}[c(e_\alpha)c(e_\beta)c(e_\gamma)c(e_\mu)c(e_\delta)c(e_\nu)]=-\sum^n_{\alpha\neq\beta,\beta\neq\mu=1}a_{1\alpha}a_{2\beta}a_{3\alpha}a_{4\mu}a_{5\beta}a_{6\mu}{\rm tr}[{\rm \texttt{id}}].\nonumber\\
\end{align}
$\mathbf{case(3-f)}$ When $\alpha\neq \beta, \alpha=\gamma, \beta\neq\mu, \beta=\nu, \mu=\delta$.\\
By (\ref{a4}), we have
\begin{align}
\label{w7}
&\sum^n_{\alpha,\beta,\gamma,\mu,\delta,\nu=1}a_{1\alpha}a_{2\beta}a_{3\gamma}a_{4\mu}a_{5\delta}a_{6\nu}{\rm tr}[c(e_\alpha)c(e_\beta)c(e_\gamma)c(e_\mu)c(e_\delta)c(e_\nu)]=\sum^n_{\alpha\neq\beta,\beta\neq\mu=1}a_{1\alpha}a_{2\beta}a_{3\alpha}a_{4\mu}a_{5\mu}a_{6\beta}{\rm tr}[{\rm \texttt{id}}].\nonumber\\
\end{align}
$\mathbf{case(3-g)}$ When $\alpha\neq \beta, \alpha=\mu, \beta=\gamma, \delta=\nu$.\\
By (\ref{a4}), we have
\begin{align}
\label{w8}
&\sum^n_{\alpha,\beta,\gamma,\mu,\delta,\nu=1}a_{1\alpha}a_{2\beta}a_{3\gamma}a_{4\mu}a_{5\delta}a_{6\nu}{\rm tr}[c(e_\alpha)c(e_\beta)c(e_\gamma)c(e_\mu)c(e_\delta)c(e_\nu)]=-\sum^n_{\alpha\neq\beta,\delta=1}a_{1\alpha}a_{2\beta}a_{3\beta}a_{4\alpha}a_{5\delta}a_{6\delta}{\rm tr}[{\rm \texttt{id}}].\nonumber\\
\end{align}
$\mathbf{case(3-h)}$ When $\alpha\neq \beta\neq\gamma, \alpha=\mu, \beta=\delta, \gamma=\nu$.\\
By (\ref{a4}), we have
\begin{align}
\label{w9}
&\sum^n_{\alpha,\beta,\gamma,\mu,\delta,\nu=1}a_{1\alpha}a_{2\beta}a_{3\gamma}a_{4\mu}a_{5\delta}a_{6\nu}{\rm tr}[c(e_\alpha)c(e_\beta)c(e_\gamma)c(e_\mu)c(e_\delta)c(e_\nu)]=\sum^n_{\alpha\neq\beta\neq\gamma=1}a_{1\alpha}a_{2\beta}a_{3\gamma}a_{4\alpha}a_{5\beta}a_{6\gamma}{\rm tr}[{\rm \texttt{id}}].\nonumber\\
\end{align}
$\mathbf{case(3-i)}$ When $\alpha\neq \beta\neq\gamma, \alpha=\mu, \beta=\nu, \delta=\gamma$.\\
By (\ref{a4}), we have
\begin{align}
\label{w10}
&\sum^n_{\alpha,\beta,\gamma,\mu,\delta,\nu=1}a_{1\alpha}a_{2\beta}a_{3\gamma}a_{4\mu}a_{5\delta}a_{6\nu}{\rm tr}[c(e_\alpha)c(e_\beta)c(e_\gamma)c(e_\mu)c(e_\delta)c(e_\nu)]=-\sum^n_{\alpha\neq\beta\neq\gamma=1}a_{1\alpha}a_{2\beta}a_{3\gamma}a_{4\alpha}a_{5\gamma}a_{6\beta}{\rm tr}[{\rm \texttt{id}}].\nonumber\\
\end{align}
$\mathbf{case(3-j)}$ When $\alpha\neq \beta, \alpha\neq\mu, \alpha=\delta, \beta=\gamma, \mu=\nu$.\\
By (\ref{a4}), we have
\begin{align}
\label{w11}
&\sum^n_{\alpha,\beta,\gamma,\mu,\delta,\nu=1}a_{1\alpha}a_{2\beta}a_{3\gamma}a_{4\mu}a_{5\delta}a_{6\nu}{\rm tr}[c(e_\alpha)c(e_\beta)c(e_\gamma)c(e_\mu)c(e_\delta)c(e_\nu)]=\sum^n_{\alpha\neq\beta,\alpha\neq\mu=1}a_{1\alpha}a_{2\beta}a_{3\beta}a_{4\mu}a_{5\alpha}a_{6\mu}{\rm tr}[{\rm \texttt{id}}].\nonumber\\
\end{align}
$\mathbf{case(3-k)}$ When $\alpha\neq \beta\neq\gamma, \alpha=\delta, \beta=\mu, \gamma=\nu$.\\
By (\ref{a4}), we have
\begin{align}
\label{w12}
&\sum^n_{\alpha,\beta,\gamma,\mu,\delta,\nu=1}a_{1\alpha}a_{2\beta}a_{3\gamma}a_{4\mu}a_{5\delta}a_{6\nu}{\rm tr}[c(e_\alpha)c(e_\beta)c(e_\gamma)c(e_\mu)c(e_\delta)c(e_\nu)]=-\sum^n_{\alpha\neq \beta\neq\gamma=1}a_{1\alpha}a_{2\beta}a_{3\gamma}a_{4\beta}a_{5\alpha}a_{6\gamma}{\rm tr}[{\rm \texttt{id}}].\nonumber\\
\end{align}
$\mathbf{case(3-l)}$ When $\alpha\neq \beta\neq\gamma, \alpha=\delta, \beta=\nu, \mu=\gamma$.\\
By (\ref{a4}), we have
\begin{align}
\label{w13}
&\sum^n_{\alpha,\beta,\gamma,\mu,\delta,\nu=1}a_{1\alpha}a_{2\beta}a_{3\gamma}a_{4\mu}a_{5\delta}a_{6\nu}{\rm tr}[c(e_\alpha)c(e_\beta)c(e_\gamma)c(e_\mu)c(e_\delta)c(e_\nu)]=\sum^n_{\alpha\neq\beta\neq\gamma=1}a_{1\alpha}a_{2\beta}a_{3\gamma}a_{4\gamma}a_{5\alpha}a_{6\beta}{\rm tr}[{\rm \texttt{id}}].\nonumber\\
\end{align}
$\mathbf{case(3-m)}$ When $\alpha\neq \beta, \alpha\neq\mu, \alpha=\nu, \beta=\gamma, \mu=\delta$.\\
By (\ref{a4}), we have
\begin{align}
\label{w14}
&\sum^n_{\alpha,\beta,\gamma,\mu,\delta,\nu=1}a_{1\alpha}a_{2\beta}a_{3\gamma}a_{4\mu}a_{5\delta}a_{6\nu}{\rm tr}[c(e_\alpha)c(e_\beta)c(e_\gamma)c(e_\mu)c(e_\delta)c(e_\nu)]=-\sum^n_{\alpha\neq \beta, \alpha\neq\mu=1}a_{1\alpha}a_{2\beta}a_{3\beta}a_{4\mu}a_{5\mu}a_{6\alpha}{\rm tr}[{\rm \texttt{id}}].\nonumber\\
\end{align}
$\mathbf{case(3-n)}$ When $\alpha\neq \beta\neq\gamma, \alpha=\nu, \beta=\mu, \gamma=\delta$.\\
By (\ref{a4}), we have
\begin{align}
\label{w15}
&\sum^n_{\alpha,\beta,\gamma,\mu,\delta,\nu=1}a_{1\alpha}a_{2\beta}a_{3\gamma}a_{4\mu}a_{5\delta}a_{6\nu}{\rm tr}[c(e_\alpha)c(e_\beta)c(e_\gamma)c(e_\mu)c(e_\delta)c(e_\nu)]=\sum^n_{\alpha\neq \beta\neq\gamma=1}a_{1\alpha}a_{2\beta}a_{3\gamma}a_{4\beta}a_{5\gamma}a_{6\alpha}{\rm tr}[{\rm \texttt{id}}].\nonumber\\
\end{align}
$\mathbf{case(3-o)}$ When $\alpha\neq \beta\neq\gamma, \alpha=\nu, \beta=\delta, \gamma=\mu$.\\
By (\ref{a4}), we have
\begin{align}
\label{w16}
&\sum^n_{\alpha,\beta,\gamma,\mu,\delta,\nu=1}a_{1\alpha}a_{2\beta}a_{3\gamma}a_{4\mu}a_{5\delta}a_{6\nu}{\rm tr}[c(e_\alpha)c(e_\beta)c(e_\gamma)c(e_\mu)c(e_\delta)c(e_\nu)]=-\sum^n_{\alpha\neq\beta\neq\gamma=1}a_{1\alpha}a_{2\beta}a_{3\gamma}a_{4\gamma}a_{5\beta}a_{6\alpha}{\rm tr}[{\rm \texttt{id}}].\nonumber\\
\end{align}
$\mathbf{case(3-p)}$ Other cases.\\
By (\ref{a4}), we have
\begin{align}
\label{w55}
&\sum^n_{\alpha,\beta,\gamma,\mu,\delta,\nu=1}a_{1\alpha}a_{2\beta}a_{3\gamma}a_{4\mu}a_{5\delta}a_{6\nu}{\rm tr}[c(e_\alpha)c(e_\beta)c(e_\gamma)c(e_\mu)c(e_\delta)c(e_\nu)]=0.\nonumber\\
\end{align}
Similar to (\ref{d45}), we have
\begin{align}
\label{w33}
&-\sum^n_{\alpha,\gamma,\delta=1}a_{1\alpha}a_{2\alpha}a_{3\gamma}a_{4\gamma}a_{5\delta}a_{6\delta}{\rm tr}[{\rm \texttt{id}}]=-g(X_1,X_2)g(X_3,X_4)g(X_5,X_6){\rm tr}[{\rm \texttt{id}}],\nonumber\\
\end{align}
\begin{align}
\label{w303}
&\sum^n_{\alpha,\gamma\neq\mu=1}a_{1\alpha}a_{2\alpha}a_{3\gamma}a_{4\mu}a_{5\gamma}a_{6\mu}{\rm tr}[{\rm \texttt{id}}]-\sum^n_{\alpha,\gamma\neq\mu=1}a_{1\alpha}a_{2\alpha}a_{3\gamma}a_{4\mu}a_{5\mu}a_{6\gamma}{\rm tr}[{\rm \texttt{id}}]\nonumber\\
&=\sum^n_{\alpha,\gamma,\mu=1}a_{1\alpha}a_{2\alpha}a_{3\gamma}a_{4\mu}a_{5\gamma}a_{6\mu}{\rm tr}[{\rm \texttt{id}}]-\sum^n_{\alpha,\gamma=\mu=1}a_{1\alpha}a_{2\alpha}a_{3\gamma}a_{4\gamma}a_{5\gamma}a_{6\gamma}{\rm tr}[{\rm \texttt{id}}]\nonumber\\
&-\sum^n_{\alpha,\gamma,\mu=1}a_{1\alpha}a_{2\alpha}a_{3\gamma}a_{4\mu}a_{5\mu}a_{6\gamma}{\rm tr}[{\rm \texttt{id}}]+\sum^n_{\alpha,\gamma=\mu=1}a_{1\alpha}a_{2\alpha}a_{3\gamma}a_{4\gamma}a_{5\gamma}a_{6\gamma}{\rm tr}[{\rm \texttt{id}}]\nonumber\\
&=g(X_1,X_2)[g(X_3,X_5)g(X_4,X_6)-g(X_3,X_6)g(X_4,X_5)]{\rm tr}[{\rm \texttt{id}}],\nonumber\\
\end{align}
\begin{align}
\label{w313}
&\sum^n_{\alpha\neq\beta,\delta=1}a_{1\alpha}a_{2\beta}a_{3\alpha}a_{4\beta}a_{5\delta}a_{6\delta}{\rm tr}[{\rm \texttt{id}}]-\sum^n_{\alpha\neq\beta,\delta=1}a_{1\alpha}a_{2\beta}a_{3\beta}a_{4\alpha}a_{5\delta}a_{6\delta}{\rm tr}[{\rm \texttt{id}}]\nonumber\\
&=[g(X_1,X_3)g(X_2,X_4)-g(X_1,X_4)g(X_2,X_3)]g(X_5,X_6){\rm tr}[{\rm \texttt{id}}],\nonumber\\
\end{align}
\begin{align}
\label{w323}
&\sum^n_{\alpha\neq\beta,\beta\neq\mu=1}a_{1\alpha}a_{2\beta}a_{3\alpha}a_{4\mu}a_{5\mu}a_{6\beta}{\rm tr}[{\rm \texttt{id}}]-\sum^n_{\alpha\neq \beta, \alpha\neq\mu=1}a_{1\alpha}a_{2\beta}a_{3\beta}a_{4\mu}a_{5\mu}a_{6\alpha}{\rm tr}[{\rm \texttt{id}}]\nonumber\\
&+\sum^n_{\alpha\neq\beta,\alpha\neq\mu=1}a_{1\alpha}a_{2\beta}a_{3\beta}a_{4\mu}a_{5\alpha}a_{6\mu}{\rm tr}[{\rm \texttt{id}}]-\sum^n_{\alpha\neq\beta,\beta\neq\mu=1}a_{1\alpha}a_{2\beta}a_{3\alpha}a_{4\mu}a_{5\beta}a_{6\mu}{\rm tr}[{\rm \texttt{id}}]\nonumber\\
&=\bigg\{g(X_1,X_3)[g(X_2,X_6)g(X_4,X_5)-g(X_2,X_5)g(X_4,X_6)]+g(X_2,X_3)[g(X_1,X_5)g(X_4,X_6)\nonumber\\
&-g(X_1,X_6)g(X_4,X_5)]\bigg\}{\rm tr}[{\rm \texttt{id}}],\nonumber\\
\end{align}
and
\begin{align}
\label{w333}
&\sum^n_{\alpha\neq\beta\neq\gamma=1}a_{1\alpha}a_{2\beta}a_{3\gamma}a_{4\alpha}a_{5\beta}a_{6\gamma}{\rm tr}[{\rm \texttt{id}}]-\sum^n_{\alpha\neq\beta\neq\gamma=1}a_{1\alpha}a_{2\beta}a_{3\gamma}a_{4\alpha}a_{5\gamma}a_{6\beta}{\rm tr}[{\rm \texttt{id}}]\nonumber\\
&-\sum^n_{\alpha\neq \beta\neq\gamma=1}a_{1\alpha}a_{2\beta}a_{3\gamma}a_{4\beta}a_{5\alpha}a_{6\gamma}{\rm tr}[{\rm \texttt{id}}]+\sum^n_{\alpha\neq\beta\neq\gamma=1}a_{1\alpha}a_{2\beta}a_{3\gamma}a_{4\gamma}a_{5\alpha}a_{6\beta}{\rm tr}[{\rm \texttt{id}}]\nonumber\\
&+\sum^n_{\alpha\neq \beta\neq\gamma=1}a_{1\alpha}a_{2\beta}a_{3\gamma}a_{4\beta}a_{5\gamma}a_{6\alpha}{\rm tr}[{\rm \texttt{id}}]-\sum^n_{\alpha\neq\beta\neq\gamma=1}a_{1\alpha}a_{2\beta}a_{3\gamma}a_{4\gamma}a_{5\beta}a_{6\alpha}{\rm tr}[{\rm \texttt{id}}]\nonumber\\
&=\bigg\{g(X_1,X_4)[g(X_2,X_5)g(X_3,X_6)-g(X_2,X_6)g(X_3,X_5)]+g(X_1,X_6)[g(X_2,X_4)g(X_3,X_5)\nonumber\\
&-g(X_2,X_5)g(X_3,X_4)]+g(X_1,X_5)[g(X_2,X_6)g(X_3,X_4)-g(X_2,X_4)g(X_3,X_6)]\bigg\}{\rm tr}[{\rm \texttt{id}}].\nonumber\\
\end{align}
Therefore
\begin{align}
\label{w66}
&{\rm tr}[c(X_1)c(X_2)c(X_3)c(X_4)c(X_5)c(X_6)]\nonumber\\
&=\bigg\{g(X_1,X_2)[g(X_3,X_5)g(X_4,X_6)-g(X_3,X_6)g(X_4,X_5)-g(X_3,X_4)g(X_5,X_6)]\nonumber\\
&+g(X_1,X_3)[g(X_2,X_4)g(X_5,X_6)-g(X_2,X_6)g(X_4,X_5)-g(X_2,X_5)g(X_3,X_6)]\nonumber\\
&+g(X_1,X_4)[g(X_2,X_5)g(X_3,X_6)-g(X_2,X_6)g(X_3,X_5)-g(X_2,X_3)g(X_5,X_6)]\nonumber\\
&+g(X_1,X_5)[g(X_2,X_6)g(X_3,X_4)-g(X_2,X_4)g(X_3,X_6)-g(X_2,X_3)g(X_4,X_6)]\nonumber\\
&+g(X_1,X_6)[g(X_2,X_4)g(X_3,X_5)-g(X_2,X_5)g(X_3,X_4)-g(X_2,X_3)g(X_4,X_5)]\bigg\}{\rm tr}[{\rm \texttt{id}}].\nonumber\\
\end{align}
When $X_j=e_j,~l=n,$ ${\rm tr}[c(e_1)\cdot\cdot\cdot c(e_n)]={\rm str}[1]=0.$ So we have
\begin{thm}\label{lko}
We have the following equalities
\begin{align}\label{jk}
&{\rm Swres}[(D^2)^{-\frac{n-2}{2}}]=0;\nonumber\\
&{\rm Wres}[c(X_1)c(X_2)(D^2)^{-\frac{n-2}{2}}]=\frac{(n-2)(4\pi)^{\frac{n}{2}}}{(\frac{n}{2}-1)!}\int_M\frac{1}{12}g(X_1,X_2)s{\rm tr}[{\rm \texttt{id}}]d{\rm Vol_{M}};\nonumber\\
&{\rm Wres}[c(X_1)c(X_2)c(X_3)c(X_4)(D^2)^{-\frac{n-2}{2}}]=\frac{(n-2)(4\pi)^{\frac{n}{2}}}{(\frac{n}{2}-1)!}\int_M-\frac{1}{12}[g(X_1,X_2)g(X_3,X_4)-g(X_1,X_3)g(X_2,X_4)\nonumber\\
&+g(X_1,X_4)g(X_2,X_3)]s{\rm tr}[{\rm \texttt{id}}]d{\rm Vol_{M}};\nonumber\\
&{\rm Wres}[c(X_1)c(X_2)c(X_3)c(X_4)c(X_5)c(X_6)(D^2)^{-\frac{n-2}{2}}]=\frac{(n-2)(4\pi)^{\frac{n}{2}}}{(\frac{n}{2}-1)!}\int_M-\frac{1}{12}\bigg\{g(X_1,X_2)[g(X_3,X_5)g(X_4,X_6)\nonumber\\
&-g(X_3,X_6)g(X_4,X_5)-g(X_3,X_4)g(X_5,X_6)]\nonumber\\
&+g(X_1,X_3)[g(X_2,X_4)g(X_5,X_6)-g(X_2,X_6)g(X_4,X_5)-g(X_2,X_5)g(X_3,X_6)]\nonumber\\
&+g(X_1,X_4)[g(X_2,X_5)g(X_3,X_6)-g(X_2,X_6)g(X_3,X_5)-g(X_2,X_3)g(X_5,X_6)]\nonumber\\
&+g(X_1,X_5)[g(X_2,X_6)g(X_3,X_4)-g(X_2,X_4)g(X_3,X_6)-g(X_2,X_3)g(X_4,X_6)]\nonumber\\
&+g(X_1,X_6)[g(X_2,X_4)g(X_3,X_5)-g(X_2,X_5)g(X_3,X_4)-g(X_2,X_3)g(X_4,X_5)]\bigg\}s{\rm tr}[{\rm \texttt{id}}]d{\rm Vol_{M}}.\nonumber\\
\end{align}
\end{thm}
\section{A Kastler-Kalau-Walze type theorem for $4$-dimensional manifolds with boundary}
\label{section:3}
 In this section, we prove the Kastler-Kalau-Walze type theorem  for the generalized noncommutative
residue $\widetilde{{\rm Wres}}[\pi^+(LD^{-1})\circ\pi^+(D^{-1})]$ on $4$-dimensional oriented compact manifolds with boundary. We firstly recall some basic facts and formulas about Boutet de
Monvel's calculus and the definition of the noncommutative residue for manifolds with boundary which will be used in the following. For more details, see Section 2 in \cite{Wa3}.\\
 \indent Let $M$ be a 4-dimensional compact oriented manifold with boundary $\partial M$.
We assume that the metric $g^{M}$ on $M$ has the following form near the boundary,
\begin{equation}
\label{b1}
g^{M}=\frac{1}{h(x_{n})}g^{\partial M}+dx _{n}^{2},
\end{equation}
where $g^{\partial M}$ is the metric on $\partial M$ and $h(x_n)\in C^{\infty}([0, 1)):=\{\widehat{h}|_{[0,1)}|\widehat{h}\in C^{\infty}((-\varepsilon,1))\}$ for
some $\varepsilon>0$ and $h(x_n)$ satisfies $h(x_n)>0$, $h(0)=1,$ where $x_n$ denotes the normal directional coordinate.\\
\indent  Then similar to \cite{Wa3}, we can compute the generalized noncommutative residue
\begin{equation}
\label{b14}
\widetilde{{\rm Wres}}[\pi^+(LD^{-1})\circ\pi^+(D^{-1})]=\int_M\int_{|\xi|=1}{\rm
trace}_{\wedge^*T^*M\bigotimes\mathbb{C}}[\sigma_{-4}(LD^{-2})]\sigma(\xi)dx+\int_{\partial M}\Phi,
\end{equation}
where
\begin{align}
\label{b15}
\Phi &=\int_{|\xi'|=1}\int^{+\infty}_{-\infty}\sum^{\infty}_{j, k=0}\sum\frac{(-i)^{|\alpha|+j+k+1}}{\alpha!(j+k+1)!}
\times {\rm trace}_{\wedge^*T^*M\bigotimes\mathbb{C}}[\partial^j_{x_n}\partial^\alpha_{\xi'}\partial^k_{\xi_n}\sigma^+_{r}(LD^{-1})(x',0,\xi',\xi_n)
\nonumber\\
&\times\partial^\alpha_{x'}\partial^{j+1}_{\xi_n}\partial^k_{x_n}\sigma_{l}(D^{-1})(x',0,\xi',\xi_n)]d\xi_n\sigma(\xi')dx',
\end{align}
and the sum is taken over $r+l-k-j-|\alpha|=-3,~~r\leq -1,~~l\leq-1$.\\

\indent By Theorem \ref{thm2}, we can compute the interior of $\widetilde{{\rm Wres}}[\pi^+(LD^{-1})\circ\pi^+(D^{-1})]$,\\
(1)when $l=2,$ we get
\begin{align}
\label{b16}
&\int_M\int_{|\xi|=1}{\rm
trace}_{\wedge^*T^*M}[\sigma_{-4}(LD^{-2})]\sigma(\xi)dx=32\pi^2\int_{M}\bigg(\frac{1}{3}g(X_1,X_2)s\bigg)d{\rm Vol_{M}}.
\end{align}
(2)when $l=4,$ we get
\begin{align}
\label{b12222}
&\int_M\int_{|\xi|=1}{\rm
trace}_{\wedge^*T^*M}[\sigma_{-4}(LD^{-2})]\sigma(\xi)dx\nonumber\\
&=32\pi^2\int_{M}\bigg(-\frac{1}{3}[g(X_1,X_2)g(X_3,X_4)-g(X_1,X_3)g(X_2,X_4)+g(X_1,X_4)g(X_2,X_3)]s\bigg)d{\rm Vol_{M}}.
\end{align}
(3)when $l=1~or~3,$ we get
\begin{align}
\label{b1202}
&\int_M\int_{|\xi|=1}{\rm
trace}_{\wedge^*T^*M}[\sigma_{-4}(LD^{-2})]\sigma(\xi)dx=0.
\end{align}
\indent Now we  need to compute $\int_{\partial M} \Phi$. By \cite{Wa3}, we get the following symbols.
\begin{lem}\label{lem2} The following identities hold:
\begin{align}
\label{b17}
&\sigma_1(D)=ic(\xi); \nonumber\\
&\sigma_0(D)=-\frac{1}{4}\sum_{i,s,t}\omega_{s,t}(e_i)c(e_i)c(e_s)c(e_t). \nonumber\\
\end{align}
\end{lem}
By $\sigma(LD^{-1})=\sigma(L)\sigma(D^{-1})$, we have
\begin{lem}\label{lem3} The following identities hold:
\begin{align}
\label{b22}
&\sigma_{-1}({D}^{-1})=\frac{ic(\xi)}{|\xi|^2};\nonumber\\
&\sigma_{-1}(L{D}^{-1})=\frac{iLc(\xi)}{|\xi|^2};\nonumber\\
&\sigma_{-2}({D}^{-1})=\frac{c(\xi)\sigma_{0}(D)c(\xi)}{|\xi|^4}+\frac{c(\xi)}{|\xi|^6}\sum_jc(dx_j)
\Big[\partial_{x_j}(c(\xi))|\xi|^2-c(\xi)\partial_{x_j}(|\xi|^2)\Big] ;\nonumber\\
&\sigma_{-2}(LD^{-1})=L\left\{\frac{c(\xi)\sigma_{0}(D)c(\xi)}{|\xi|^4}+\frac{c(\xi)}{|\xi|^6}\sum_jc(dx_j)
\Big[\partial_{x_j}(c(\xi))|\xi|^2-c(\xi)\partial_{x_j}(|\xi|^2)\Big]\right\}.
\end{align}
\end{lem}
\indent When $n=4$, then ${\rm tr}_{S(TM)}[{\rm \texttt{id}}]={\rm dim}(\wedge^*(\mathbb{R}^2))=4$, the sum is taken over $
r+\ell-k-j-|\alpha|=-3,~~r\leq -1,~~\ell\leq-1,$ then we have the following five cases:
~\\
\noindent  {\bf case a)~I)}~$r=-1,~\ell=-1,~k=j=0,~|\alpha|=1$.\\
\noindent By (\ref{b15}), we get
\begin{equation}
\label{b24}
\Phi_1=-\int_{|\xi'|=1}\int^{+\infty}_{-\infty}\sum_{|\alpha|=1}
 {\rm tr}[\partial^\alpha_{\xi'}\pi^+_{\xi_n}\sigma_{-1}(L{D}^{-1})\times
 \partial^\alpha_{x'}\partial_{\xi_n}\sigma_{-1}(D^{-1})](x_0)d\xi_n\sigma(\xi')dx'.
\end{equation}
By Lemma 2.2 in \cite{Wa3}, for $i<n$, then
\begin{equation}
\label{b25}\partial_{x_i}\left(\frac{ic(\xi)}{|\xi|^2}\right)(x_0)=
\frac{i\partial_{x_i}[c(\xi)](x_0)}{|\xi|^2}
-\frac{ic(\xi)\partial_{x_i}(|\xi|^2)(x_0)}{|\xi|^4}=0,
\end{equation}
\noindent so when $l$ is an integer between 1 and 4, we have $\Phi_1=0$.\\
 \noindent  {\bf case a)~II)}~$r=-1,~\ell=-1,~k=|\alpha|=0,~j=1$.\\
\noindent By (\ref{b15}), we get
\begin{equation}
\label{b26}
\Phi_2=-\frac{1}{2}\int_{|\xi'|=1}\int^{+\infty}_{-\infty} {\rm
trace} [\partial_{x_n}\pi^+_{\xi_n}\sigma_{-1}(LD^{-1})\times
\partial_{\xi_n}^2\sigma_{-1}(D^{-1})](x_0)d\xi_n\sigma(\xi')dx'.
\end{equation}
\noindent By Lemma \ref{lem3}, we have\\
\begin{eqnarray}\label{b27}\partial^2_{\xi_n}\sigma_{-1}(D^{-1})(x_0)=i\left(-\frac{6\xi_nc(dx_n)+2c(\xi')}
{|\xi|^4}+\frac{8\xi_n^2c(\xi)}{|\xi|^6}\right);
\end{eqnarray}
When $l=2$, we have
\begin{align}\label{b28}
\partial_{x_n}\sigma_{-1}(LD^{-1})(x_0)&=\partial_{x_n}\sigma_{-1}(c(X_1)c(X_2)D^{-1})(x_0)\nonumber\\
&=\frac{i\sum^4_{\alpha,\beta=1}\partial_{x_n}(a_{1\alpha})a_{2\beta}c(e_\alpha)c(e_\beta)c(\xi)}{|\xi|^2}+\frac{i\sum^4_{\alpha,\beta=1}a_{1\alpha}\partial_{x_n}(a_{2\beta})c(e_\alpha)c(e_\beta)c(\xi)}{|\xi|^2}\nonumber\\
&+\frac{i\sum^4_{\alpha,\beta=1}a_{1\alpha}a_{2\beta}c(e_\alpha)c(e_\beta)\partial_{x_n}c(\xi')(x_0)}{|\xi|^2}-\frac{i\sum^4_{\alpha,\beta=1}a_{1\alpha}a_{2\beta}c(e_\alpha)c(e_\beta)c(\xi)|\xi'|^2h'(0)}{|\xi|^4}.
\end{align}
By (2.1.1), (2.1.2) in \cite{Wa3} and the Cauchy integral formula, we have
\begin{align}\label{b29}
\pi^+_{\xi_n}\left[\frac{i\sum^4_{\alpha,\beta=1}\partial_{x_n}(a_{1\alpha})a_{2\beta}c(e_\alpha)c(e_\beta)c(\xi)}{|\xi|^2}\right](x_0)|_{|\xi'|=1}&=i\sum^4_{\alpha,\beta=1}\partial_{x_n}(a_{1\alpha})a_{2\beta}c(e_\alpha)c(e_\beta)\pi^+_{\xi_n}\left[\frac{c(\xi)}{|\xi|^2}\right](x_0)|_{|\xi'|=1}\nonumber\\
&=\sum^4_{\alpha,\beta=1}\partial_{x_n}(a_{1\alpha})a_{2\beta}c(e_\alpha)c(e_\beta)\frac{c(\xi')+ic(dx_n)}{2(\xi_n-i)}.\nonumber\\
\end{align}
Similarly, we have,
\begin{align}\label{b209}
\pi^+_{\xi_n}\left[\frac{i\sum^4_{\alpha,\beta=1}a_{1\alpha}\partial_{x_n}(a_{2\beta})c(e_\alpha)c(e_\beta)c(\xi)}{|\xi|^2}\right](x_0)|_{|\xi'|=1}&=\frac{\sum^4_{\alpha,\beta=1}a_{1\alpha}\partial_{x_n}(a_{2\beta})c(e_\alpha)c(e_\beta)[c(\xi')+ic(dx_n)]}{2(\xi_n-i)}.\nonumber\\
\end{align}
\begin{align}\label{30}
\pi^+_{\xi_n}\left[\frac{i\sum^4_{\alpha,\beta=1}a_{1\alpha}a_{2\beta}c(e_\alpha)c(e_\beta)\partial_{x_n}c(\xi')}{|\xi|^2}\right](x_0)|_{|\xi'|=1}=\frac{\sum^4_{\alpha,\beta=1}a_{1\alpha}a_{2\beta}c(e_\alpha)c(e_\beta)\partial_{x_n}[c(\xi')](x_0)}{2(\xi_n-i)}.\nonumber\\
\end{align}
\begin{align}\label{300}
\pi^+_{\xi_n}\left[\frac{i\sum^4_{\alpha,\beta=1}a_{1\alpha}a_{2\beta}c(e_\alpha)c(e_\beta)c(\xi)|\xi'|^2h'(0)}{|\xi|^4}\right](x_0)|_{|\xi'|=1}=-i\sum^4_{\alpha,\beta=1}a_{1\alpha}a_{2\beta}c(e_\alpha)c(e_\beta)\left[\frac{(i\xi_n+2)c(\xi')+ic(dx_n)}{4(\xi_n-i)^2}\right].\nonumber\\
\end{align}
\noindent By the relation of the Clifford action and ${\rm tr}{ab}={\rm tr }{ba}$, we have the equalities:\\
\begin{align}\label{32}
&{\rm tr}[c(\xi')c(dx_n)]=0;~~{\rm tr}[c(dx_n)^2]=-4;~~{\rm tr}[c(\xi')^2](x_0)|_{|\xi'|=1}=-4;~~{\rm tr}[\partial_{x_n}c(\xi')c(dx_n)]=0;\nonumber\\
&{\rm tr}[\partial_{x_n}c(\xi')c(\xi')](x_0)|_{|\xi'|=1}=-2h'(0);~~\sum_{\alpha,\beta=1}^4\partial_{x_n}(a_{1\alpha})(a_{2\alpha})+\sum_{\alpha,\beta=1}^4a_{1\alpha}\partial_{x_n}(a_{2\alpha})=\partial_{x_n}[g(X_1,X_2)];\nonumber\\
&{\rm tr}[\sum_{\alpha,\beta=1}^4a_{1\alpha}a_{2\beta}c(e_\alpha)c(e_\beta)c(\xi')c(\xi')]={\rm tr}[\sum_{\alpha,\beta=1}^4a_{1\alpha}a_{2\beta}c(e_\alpha)c(e_\beta)c(dx_n)c(dx_n)]=4g(X_1,X_2);\nonumber\\
&{\rm tr}[\sum_{\alpha,\beta=1}^4\partial_{x_n}(a_{1\alpha})a_{2\beta}c(e_\alpha)c(e_\beta)c(\xi')c(\xi')]={\rm tr}[\sum_{\alpha,\beta=1}^4\partial_{x_n}(a_{1\alpha})a_{2\beta}c(e_\alpha)c(e_\beta)c(dx_n)c(dx_n)]=4\sum_{\alpha,\beta=1}^4\partial_{x_n}(a_{1\alpha})a_{2\alpha}.\nonumber\\
\end{align}
We note that $i<n,~\int_{|\xi'|=1}\{\xi_{i_{1}}\xi_{i_{2}}\cdots\xi_{i_{2d+1}}\}\sigma(\xi')=0$,
so ${\rm tr }[c(e_\alpha)c(e_\beta)c(\xi')c(dx_n)]$ has no contribution for computing {\bf case a)~II)}.\\
Then, we have
\begin{align}\label{33}
&{\rm
trace} [\partial_{x_n}\pi^+_{\xi_n}\sigma_{-1}(c(X_1)c(X_2)D^{-1})\times
\partial_{\xi_n}^2\sigma_{-1}(D^{-1})](x_0)\nonumber\\
&=4\frac{3\xi_n^2i-\xi_n^3-3\xi_n+i}{(\xi_n-i)^4(\xi_n+i)^3}\sum^4_{\alpha=1}\partial_{x_n}[g(X_1,X_2)]+2h'(0)\frac{8\xi_n^3i+5\xi_ni+3+11\xi_n^2-\xi_n^3}{(\xi_n-i)^5(\xi_n+i)^3}g(X_1,X_2).\nonumber\\
\end{align}
Then\\
\begin{align}\label{35}
\Phi_2&=-\frac{1}{2}\int_{|\xi'|=1}\int^{+\infty}_{-\infty}\bigg\{4\frac{3\xi_n^2i-\xi_n^3-3\xi_n+i}{(\xi_n-i)^4(\xi_n+i)^3}\partial_{x_n}[g(X_1,X_2)]+2\frac{8\xi_n^3i+5\xi_ni+3+11\xi_n^2-\xi_n^3}{(\xi_n-i)^5(\xi_n+i)^3}g(X_1,X_2)\bigg\}d\xi_n\sigma(\xi')dx'\nonumber\\
&=-2\partial_{x_n}[g(X_1,X_2)]\Omega_3\int_{\Gamma^{+}}\frac{3\xi_n^2i-\xi_n^3-3\xi_n+i}{(\xi_n-i)^4(\xi_n+i)^3}d\xi_{n}dx'-g(X_1,X_2)\Omega_3\int_{\Gamma^{+}}\frac{8\xi_n^3i+5\xi_ni+3+11\xi_n^2-\xi_n^3}{(\xi_n-i)^5(\xi_n+i)^3}d\xi_{n}dx'\nonumber\\
&=-2\partial_{x_n}[g(X_1,X_2)]\Omega_3\frac{2\pi i}{3!}\left[\frac{3\xi_n^2i-\xi_n^3-3\xi_n+i}{(\xi_n+i)^3}\right]^{(3)}\bigg|_{\xi_n=i}dx'\nonumber\\
&-g(X_1,X_2)\Omega_3\frac{2\pi i}{4!}\left[\frac{8\xi_n^3i+5\xi_ni+3+11\xi_n^2-\xi_n^3}{(\xi_n+i)^3}\right]^{(4)}\bigg|_{\xi_n=i}dx'\nonumber\\
&=\left\{-\frac{\partial_{x_n}[g(X_1,X_2)]}{2}-\frac{3h'(0)g(X_1,X_2)}{8}\right\}\pi\Omega_3dx',\nonumber\\
\end{align}
where ${\rm \Omega_{3}}$ is the canonical volume of $S^{3}.$\\
Similar to (\ref{a286})-(\ref{d45}), we get\\
\begin{align}\label{372}
&{\rm tr}[\sum_{\alpha,\beta,\gamma,\mu=1}^4\partial_{x_n}(a_{1\alpha})a_{2\beta}a_{3\gamma}a_{4\mu}c(e_\alpha)c(e_\beta)c(e_\gamma)c(e_\mu)c(\xi')c(\xi')]\nonumber\\
&={\rm tr}[\sum_{\alpha,\beta,\gamma,\mu=1}^4\partial_{x_n}(a_{1\alpha})a_{2\beta}a_{3\gamma}a_{4\mu}c(e_\alpha)c(e_\beta)c(e_\gamma)c(e_\mu)c(dx_n)c(dx_n)]\nonumber\\
&=4(\sum^n_{\alpha,\beta=1}\partial_{x_n}(a_{1\alpha})a_{2\beta}a_{3\alpha}a_{4\beta}-\sum^n_{\alpha,\gamma=1}\partial_{x_n}(a_{1\alpha})a_{2\alpha}a_{3\gamma}a_{4\gamma}-\sum^n_{\alpha,\beta=1}\partial_{x_n}(a_{1\alpha})a_{2\beta}a_{3\beta}a_{4\alpha}).\nonumber\\
\end{align}
Set
\begin{align}\label{382}
&A^1=\sum^n_{\alpha,\beta=1}\partial_{x_n}(a_{1\alpha})a_{2\beta}a_{3\alpha}a_{4\beta}-\sum^n_{\alpha,\gamma=1}\partial_{x_n}(a_{1\alpha})a_{2\alpha}a_{3\gamma}a_{4\gamma}-\sum^n_{\alpha,\beta=1}\partial_{x_n}(a_{1\alpha})a_{2\beta}a_{3\beta}a_{4\alpha};\nonumber\\
&A^2=\sum^n_{\alpha,\beta=1}a_{1\alpha}\partial_{x_n}(a_{2\beta})a_{3\alpha}a_{4\beta}-\sum^n_{\alpha,\gamma=1}a_{1\alpha}\partial_{x_n}(a_{2\alpha})a_{3\gamma}a_{4\gamma}-\sum^n_{\alpha,\beta=1}a_{1\alpha}\partial_{x_n}(a_{2\beta})a_{3\beta}a_{4\alpha};\nonumber\\
&A^3=\sum^n_{\alpha,\beta=1}a_{1\alpha}a_{2\beta}\partial_{x_n}(a_{3\alpha})a_{4\beta}-\sum^n_{\alpha,\gamma=1}a_{1\alpha}a_{2\alpha}\partial_{x_n}(a_{3\gamma})a_{4\gamma}-\sum^n_{\alpha,\beta=1}a_{1\alpha}a_{2\beta}\partial_{x_n}(a_{3\beta})a_{4\alpha};\nonumber\\
&A^4=\sum^n_{\alpha,\beta=1}a_{1\alpha}a_{2\beta}a_{3\alpha}\partial_{x_n}(a_{4\beta})-\sum^n_{\alpha,\gamma=1}a_{1\alpha}a_{2\alpha}a_{3\gamma}\partial_{x_n}(a_{4\gamma})-\sum^n_{\alpha,\beta=1}a_{1\alpha}a_{2\beta}a_{3\beta}\partial_{x_n}(a_{4\alpha}).\nonumber\\
\end{align}
By
\begin{align}\label{3002}
&\sum^n_{\alpha,\gamma=1}\partial_{x_n}(a_{1\alpha})a_{2\alpha}a_{3\gamma}a_{4\gamma}+\sum^n_{\alpha,\gamma=1}a_{1\alpha}\partial_{x_n}(a_{2\alpha})a_{3\gamma}a_{4\gamma}+\sum^n_{\alpha,\gamma=1}a_{1\alpha}a_{2\alpha}\partial_{x_n}(a_{3\gamma})a_{4\gamma}+\sum^n_{\alpha,\gamma=1}a_{1\alpha}a_{2\alpha}a_{3\gamma}\partial_{x_n}(a_{4\gamma})\nonumber\\
&=\sum^n_{\alpha,\gamma=1}\partial_{x_n}(a_{1\alpha}a_{2\alpha})a_{3\gamma}a_{4\gamma}+\sum^n_{\alpha,\gamma=1}a_{1\alpha}a_{2\alpha}\partial_{x_n}(a_{3\gamma}a_{4\gamma})\nonumber\\
&=\partial_{x_n}[g(X_1,X_2)g(X_3,X_4)];\nonumber\\
&\sum^n_{\alpha,\beta=1}\partial_{x_n}(a_{1\alpha})a_{2\beta}a_{3\alpha}a_{4\beta}+\sum^n_{\alpha,\beta=1}a_{1\alpha}\partial_{x_n}(a_{2\beta})a_{3\alpha}a_{4\beta}+\sum^n_{\alpha,\beta=1}a_{1\alpha}a_{2\beta}\partial_{x_n}(a_{3\alpha})a_{4\beta}+\sum^n_{\alpha,\beta=1}a_{1\alpha}a_{2\beta}a_{3\alpha}\partial_{x_n}(a_{4\beta})\nonumber\\
&=\sum^n_{\alpha,\beta=1}\partial_{x_n}(a_{1\alpha}a_{3\alpha})a_{2\beta}a_{4\beta}+\sum^n_{\alpha,\beta=1}a_{1\alpha}a_{3\alpha}\partial_{x_n}(a_{2\beta}a_{4\beta})\nonumber\\
&=\partial_{x_n}[g(X_1,X_3)g(X_2,X_4)];\nonumber\\
&\sum^n_{\alpha,\beta=1}\partial_{x_n}(a_{1\alpha})a_{2\beta}a_{3\beta}a_{4\alpha}+\sum^n_{\alpha,\beta=1}a_{1\alpha}\partial_{x_n}(a_{2\beta})a_{3\beta}a_{4\alpha}+\sum^n_{\alpha,\beta=1}a_{1\alpha}a_{2\beta}\partial_{x_n}(a_{3\beta})a_{4\alpha}+\sum^n_{\alpha,\beta=1}a_{1\alpha}a_{2\beta}a_{3\beta}\partial_{x_n}(a_{4\alpha})\nonumber\\
&=\sum^n_{\alpha,\beta=1}\partial_{x_n}(a_{1\alpha}a_{4\alpha})a_{2\beta}a_{3\beta}+\sum^n_{\alpha,\beta=1}a_{1\alpha}a_{4\alpha}\partial_{x_n}(a_{2\beta}a_{3\beta})\nonumber\\
&=\partial_{x_n}[g(X_1,X_4)g(X_2,X_3)].\nonumber\\
\end{align}
Then, we have
\begin{align}\label{3225}
A^1+A^2+A^3+A^4=\partial_{x_n}[g(X_1,X_3)g(X_2,X_4)-g(X_1,X_4)g(X_2,X_3)-g(X_1,X_2)g(X_3,X_4)].\nonumber\\
\end{align}
So when $l=4,$ we get
\begin{align}\label{uuu}
\Phi_2&=\bigg\{-\frac{1}{2}\partial_{x_n}[g(X_1,X_3)g(X_2,X_4)-g(X_1,X_4)g(X_2,X_3)-g(X_1,X_2)g(X_3,X_4)]+\frac{3h'(0)}{8}\nonumber\\
&[g(X_1,X_4)g(X_2,X_3)-g(X_1,X_2)g(X_3,X_4)+g(X_1,X_3)g(X_2,X_4)]\bigg\}\pi\Omega_3dx'.\nonumber\\
\end{align}
\noindent  {\bf case a)~III)}~$r=-1,~\ell=-1,~j=|\alpha|=0,~k=1$.\\
\noindent By (\ref{b15}), we get
\begin{equation}\label{36}
\Phi_3=-\frac{1}{2}\int_{|\xi'|=1}\int^{+\infty}_{-\infty}
{\rm trace} [\partial_{\xi_n}\pi^+_{\xi_n}\sigma_{-1}(LD^{-1})\times
\partial_{\xi_n}\partial_{x_n}\sigma_{-1}(D^{-1})](x_0)d\xi_n\sigma(\xi')dx'.
\end{equation}
\noindent By Lemma \ref{lem3}, we have\\
\begin{eqnarray}\label{37}
\partial_{\xi_n}\partial_{x_n}\sigma_{-1}(D^{-1})(x_0)|_{|\xi'|=1}
=-ih'(0)\left[\frac{c(dx_n)}{|\xi|^4}-4\xi_n\frac{c(\xi')
+\xi_nc(dx_n)}{|\xi|^6}\right]-\frac{2\xi_ni\partial_{x_n}[c(\xi')](x_0)}{|\xi|^4};
\end{eqnarray}
\begin{eqnarray}\label{38}
\partial_{\xi_n}\pi^+_{\xi_n}\sigma_{-1}(LD^{-1})(x_0)|_{|\xi'|=1}=-\frac{\sum^4_{\alpha,\beta=1}a_{1\alpha}a_{2\beta}c(e_\alpha)c(e_\beta)[c(\xi')+ic(dx_n)]}{2(\xi_n-i)^2}.
\end{eqnarray}
Similar to {\rm case~a)~II)}, when $l=2,$ we have\\
\begin{align}\label{39}
&{\rm trace} [\partial_{\xi_n}\pi^+_{\xi_n}\sigma_{-1}(LD^{-1})\times
\partial_{\xi_n}\partial_{x_n}\sigma_{-1}(D^{-1})](x_0)=2h'(0)\frac{-5i\xi_n+3\xi_n^2+\xi_n^3i+1}{(\xi_n-i)^5(\xi_n+i)^3}g(X_1,X_2).\nonumber\\
\end{align}
So we have
\begin{align}\label{41}
\Phi_3&=-\int_{|\xi'|=1}\int^{+\infty}_{-\infty}h'(0)\frac{-5i\xi_n+3\xi_n^2+\xi_n^3i+1}{(\xi_n-i)^5(\xi_n+i)^3}g(X_1,X_2)d\xi_n\sigma(\xi')dx'\nonumber\\
&=-h'(0)g(X_1,X_2)\Omega_3\int_{\Gamma^{+}}\frac{-5i\xi_n+3\xi_n^2+\xi_n^3i+1}{(\xi_n-i)^5(\xi_n+i)^3}d\xi_{n}dx'\nonumber\\
&=-h'(0)g(X_1,X_2)\Omega_3\frac{2\pi i}{4!}\left[\frac{-5i\xi_n+3\xi_n^2+\xi_n^3i+1}{(\xi_n+i)^3}\right]^{(4)}\bigg|_{\xi_n=i}dx'\nonumber\\
&=-\frac{3h'(0)}{8}g(X_1,X_2)\pi\Omega_3dx'.\nonumber\\
\end{align}
Similarly, when $l=4,$ we get
\begin{align}\label{1uu}
\Phi_3&=\frac{3h'(0)}{8}[g(X_1,X_2)g(X_3,X_4)-g(X_1,X_3)g(X_2,X_4)+g(X_1,X_4)g(X_2,X_3)]\pi\Omega_3dx'.\nonumber\\
\end{align}
\noindent  {\bf case b)}~$r=-2,~\ell=-1,~k=j=|\alpha|=0$.\\
\noindent By (\ref{b15}), we get
\begin{align}\label{42}
\Phi_4&=-i\int_{|\xi'|=1}\int^{+\infty}_{-\infty}{\rm trace} [\pi^+_{\xi_n}\sigma_{-2}(LD^{-1})\times
\partial_{\xi_n}\sigma_{-1}(D^{-1})](x_0)d\xi_n\sigma(\xi')dx'.
\end{align}
 By Lemma \ref{lem3}, we have\\
\begin{align}\label{43}
\sigma_{-2}(LD^{-1})(x_0)=L\left\{\frac{c(\xi)\sigma_{0}(D)c(\xi)}{|\xi|^4}+\frac{c(\xi)}{|\xi|^6}\sum_jc(dx_j)
\Big[\partial_{x_j}(c(\xi))|\xi|^2-c(\xi)\partial_{x_j}(|\xi|^2)\Big]\right\},
\end{align}
where
\begin{align}\label{44}
\sigma_{0}(D)(x_0)&=-\frac{1}{4}\sum_{s,t,i}\omega_{s,t}(e_i)
(x_{0})c(e_i)c(e_s)c(e_t).\nonumber\\
\end{align}
We denote
\begin{align}\label{45}
H(x_0)&=-\frac{1}{4}\sum_{s,t,i}\omega_{s,t}(e_i)
(x_{0})c(e_i)c(e_s)c(e_t),
\end{align}
where $H(x_0)=c_0c(dx_n)$ and $c_0=-\frac{3}{4}h'(0)$.\\
Then
\begin{align}\label{46}
\pi^+_{\xi_n}\sigma_{-2}(LD^{-1}(x_0))|_{|\xi'|=1}&=
\pi^+_{\xi_n}\Big[\frac{Lc(\xi)H(x_0)c(\xi)+Lc(\xi)c(dx_n)\partial_{x_n}[c(\xi')](x_0)}{(1+\xi_n^2)^2}-h'(0)\frac{Lc(\xi)c(dx_n)c(\xi)}{(1+\xi_n^{2})^3}\Big].
\end{align}
Since
\begin{align}\label{50}
\partial_{\xi_n}\sigma_{-1}(D^{-1})=i\left[\frac{c(dx_n)}{1+\xi_n^2}-\frac{2\xi_nc(\xi')+2\xi_n^2c(dx_n)}{(1+\xi_n^2)^2}\right].
\end{align}
By computations, we have
\begin{eqnarray}\label{52}
\pi^+_{\xi_n}\Big[\frac{Lc(\xi)H(x_0)c(\xi)+Lc(\xi)c(dx_n)\partial_{x_n}[c(\xi')](x_0)}{(1+\xi_n^2)^2}\Big]-h'(0)\pi^+_{\xi_n}\Big[L\frac{c(\xi)c(dx_n)c(\xi)}{(1+\xi_n)^3}\Big]:= E_1-E_2,
\end{eqnarray}
where
\begin{align}\label{53}
E_1&=\frac{-1}{4(\xi_n-i)^2}[(2+i\xi_n)Lc(\xi')H(x_0)c(\xi')+i\xi_nLc(dx_n)H(x_0)c(dx_n)\nonumber\\
&+(2+i\xi_n)Lc(\xi')c(dx_n)\partial_{x_n}c(\xi')+iLc(dx_n)H(x_0)c(\xi')
+iLc(\xi')H(x_0)c(dx_n)-iL\partial_{x_n}c(\xi')]
\end{align}
and
\begin{align}\label{54}
E_2&=\frac{h'(0)}{2}L\left[\frac{c(dx_n)}{4i(\xi_n-i)}+\frac{c(dx_n)-ic(\xi')}{8(\xi_n-i)^2}
+\frac{3\xi_n-7i}{8(\xi_n-i)^3}[ic(\xi')-c(dx_n)]\right].
\end{align}
By (\ref{50}) and (\ref{54}), when $l=2,$ we have\\
\begin{eqnarray}\label{55}{\rm tr }[E_2\times\partial_{\xi_n}\sigma_{-1}(D^{-1})]|_{|\xi'|=1}
&=-2ih'(0)\frac{-i\xi_n^2-\xi_n+4i}{4(\xi_n-i)^3(\xi_n+i)^2}g(X_1,X_2).
\end{eqnarray}
Then,
\begin{align}\label{567}
&-i\int_{|\xi'|=1}\int^{+\infty}_{-\infty}{\rm trace} [E_2\times
\partial_{\xi_n}\sigma_{-1}(D^{-1})](x_0)d\xi_n\sigma(\xi')dx'\nonumber\\
&=-i\int_{|\xi'|=1}\int^{+\infty}_{-\infty}-2ih'(0)\frac{-i\xi_n^2-\xi_n+4i}{4(\xi_n-i)^3(\xi_n+i)^2}g(X_1,X_2)d\xi_n\sigma(\xi')dx'\nonumber\\
&=-\frac{1}{2}h'(0)g(X_1,X_2)\Omega_3\int_{\Gamma^{+}}\frac{-i\xi_n^2-\xi_n+4i}{4(\xi_n-i)^3(\xi_n+i)^2}d\xi_{n}dx'\nonumber\\
&=-\frac{1}{2}h'(0)g(X_1,X_2)\Omega_3\frac{2\pi i}{3!}\left[\frac{-i\xi_n^2-\xi_n+4i}{(\xi_n+i)^2}\right]^{(3)}\bigg|_{\xi_n=i}dx'\nonumber\\
&=\frac{3h'(0)}{4}g(X_1,X_2)\pi\Omega_3dx'.\nonumber\\
\end{align}
By (\ref{50}) and (\ref{53}), when $l=2,$ we have
\begin{align}\label{56}
{\rm tr }[E_1\times\partial_{\xi_n}\sigma_{-1}(D^{-1})]|_{|\xi'|=1}=
\frac{-3h'(0)i}{2(\xi_n-i)^2(\xi_n+i)^2}g(X_1,X_2)-\frac{(\xi_n^2-i\xi_n-2)h'(0)}{2(\xi_n-i)^3(\xi_n+i)^2}g(X_1,X_2).
\end{align}
By (\ref{56}), we have
\begin{align}\label{57}
&-i\int_{|\xi'|=1}\int^{+\infty}_{-\infty}{\rm trace} [E_1\times
\partial_{\xi_n}\sigma_{-1}(D^{-1})](x_0)d\xi_n\sigma(\xi')dx'\nonumber\\
&=-i\int_{|\xi'|=1}\int^{+\infty}_{-\infty}\frac{-3h'(0)i}{2(\xi_n-i)^2(\xi_n+i)^2}g(X_1,X_2)d\xi_n\sigma(\xi')dx'\nonumber\\
&-i\int_{|\xi'|=1}\int^{+\infty}_{-\infty}\frac{(\xi_n^2-i\xi_n-2)h'(0)}{2(\xi_n-i)^3(\xi_n+i)^2}g(X_1,X_2)d\xi_n\sigma(\xi')dx'\nonumber\\
&=-\frac{3}{2}h'(0)g(X_1,X_2)\Omega_3\int_{\Gamma^{+}}\frac{1}{(\xi_n-i)^2(\xi_n+i)^2}d\xi_{n}dx'-\frac{i}{2}h'(0)g(X_1,X_2)\Omega_3\int_{\Gamma^{+}}\frac{\xi_n^2-i\xi_n-2}{(\xi_n-i)^3(\xi_n+i)^2}d\xi_{n}dx'\nonumber\\
&=-\frac{3}{2}h'(0)g(X_1,X_2)\Omega_3\frac{2\pi i}{1!}\left[\frac{1}{(\xi_n+i)^2}\right]^{(1)}\bigg|_{\xi_n=i}dx'-\frac{i}{2}h'(0)g(X_1,X_2)\Omega_3\frac{2\pi i}{3!}\left[\frac{\xi_n^2-i\xi_n-2}{(\xi_n+i)^2}\right]^{(2)}\bigg|_{\xi_n=i}dx'\nonumber\\
&=-\frac{15h'(0)}{8}g(X_1,X_2)\pi\Omega_3dx'.\nonumber\\
\end{align}
Then,  when $l=2,$ we have\\
\begin{align}\label{60}
\Phi_4=-\frac{15h'(0)}{8}g(X_1,X_2)\pi\Omega_3dx'+\frac{3h'(0)}{4}g(X_1,X_2)\pi\Omega_3dx'=-\frac{9h'(0)}{8}g(X_1,X_2)\pi \Omega_3dx'.
\end{align}
Similarly, when $l=4,$ we get
\begin{align}\label{15u}
\Phi_4&=\frac{9h'(0)}{8}[g(X_1,X_2)g(X_3,X_4)-g(X_1,X_3)g(X_2,X_4)+g(X_1,X_4)g(X_2,X_3)]\pi \Omega_3dx'.\nonumber\\
\end{align}
\noindent {\bf  case c)}~$r=-1,~\ell=-2,~k=j=|\alpha|=0$.\\
By (\ref{b15}), we get
\begin{align}\label{61}
\Phi_5=-i\int_{|\xi'|=1}\int^{+\infty}_{-\infty}{\rm trace} [\pi^+_{\xi_n}\sigma_{-1}(LD^{-1})\times
\partial_{\xi_n}\sigma_{-2}(D^{-1})](x_0)d\xi_n\sigma(\xi')dx'.
\end{align}
By Lemma \ref{lem3}, we have
\begin{align}\label{62}
\pi^+_{\xi_n}\sigma_{-1}(LD^{-1})|_{|\xi'|=1}=\frac{L[c(\xi')+ic(dx_n)]}{2(\xi_n-i)}.
\end{align}
Since
\begin{equation}\label{63}
\sigma_{-2}(D^{-1})(x_0)=\frac{c(\xi)\sigma_{0}(D)(x_0)c(\xi)}{|\xi|^4}+\frac{c(\xi)}{|\xi|^6}c(dx_n)
\bigg[\partial_{x_n}[c(\xi')](x_0)|\xi|^2-c(\xi)h'(0)|\xi|^2_{\partial_
M}\bigg],
\end{equation}
where
\begin{align}\label{64}
\sigma_{0}(D)(x_0)&=-\frac{1}{4}\sum_{s,t,i}\omega_{s,t}(e_i)(x_{0})c(e_i)c(e_s)c(e_t).
\end{align}
By computations, we have
\begin{align}\label{68}
&\partial_{\xi_n}\sigma_{-2}(D^{-1})(x_0)|_{|\xi'|=1}=\partial_{\xi_n}\bigg\{\frac{c(\xi)H(x_0)c(\xi)}{|\xi|^4}+\frac{c(\xi)}{|\xi|^6}c(dx_n)[\partial_{x_n}[c(\xi')](x_0)|\xi|^2-c(\xi)h'(0)]\bigg\}\nonumber\\
&=\frac{1}{(1+\xi_n^2)^3}\bigg[(2\xi_n-2\xi_n^3)c(dx_n)Hc(dx_n)
+(1-3\xi_n^2)c(dx_n)Hc(\xi')\nonumber\\
&+(1-3\xi_n^2)c(\xi')Hc(dx_n)
-4\xi_nc(\xi')Hc(\xi')
+(3\xi_n^2-1)\partial_{x_n}c(\xi')\nonumber\\
&-4\xi_nc(\xi')c(dx_n){\partial}_{x_n}c(\xi')
+2h'(0)c(\xi')+2h'(0)\xi_nc(dx_n)\bigg]+6\xi_nh'(0)\frac{c(\xi)c(dx_n)c(\xi)}{(1+\xi^2_n)^4}.
\end{align}
By (\ref{62}) and (\ref{68}), we have
\begin{align}\label{71}
{\rm tr}[\pi^+_{\xi_n}\sigma_{-1}(LD^{-1})\times
\partial_{\xi_n}\sigma_{-2}(D^{-1})](x_0)|_{|\xi'|=1}
=-\frac{12h'(0)i\xi_n}{(\xi_n-i)^3(\xi_n+i)^4}g(X_1,X_2)-\frac{3h'(0)(i\xi_n^2+\xi_n-2i)}{(\xi_n-i)^3(\xi_n+i)^3}g(X_1,X_2).\nonumber\\
\end{align}
So, we have
\begin{align}\label{74}
\Phi_5&=-i\int_{|\xi'|=1}\int^{+\infty}_{-\infty}{\rm tr}[\pi^+_{\xi_n}\sigma_{-1}(LD^{-1})\times
\partial_{\xi_n}\sigma_{-2}(D^{-1})](x_0)d\xi_n\sigma(\xi')dx'\nonumber\\
&=-i\int_{|\xi'|=1}\int^{+\infty}_{-\infty}-\frac{12h'(0)i\xi_n}{(\xi_n-i)^3(\xi_n+i)^4}g(X_1,X_2)d\xi_n\sigma(\xi')dx'\nonumber\\
&+i\int_{|\xi'|=1}\int^{+\infty}_{-\infty}\frac{3h'(0)(i\xi_n^2+\xi_n-2i)}{(\xi_n-i)^3(\xi_n+i)^3}g(X_1,X_2)d\xi_n\sigma(\xi')dx'\nonumber\\
&=-12g(X_1,X_2)\Omega_3\int_{\Gamma^{+}}\frac{\xi_n}{(\xi-i)^3(\xi+i)^4}d\xi_{n}dx'+3ih'(0)g(X_1,X_2)\Omega_3\int_{\Gamma^{+}}\frac{i\xi_n^2+\xi_n-2i}{(\xi-i)^3(\xi+i)^3}d\xi_{n}dx'\nonumber\\
&=-\frac{ih'(0)}{2}g(X_1,X_2)\Omega_3\frac{2\pi i}{2!}\left[\frac{\xi_n}{(\xi+i)^4}\right]^{(2)}\bigg|_{\xi_n=i}dx'+3ih'(0)g(X_1,X_2)\Omega_3\frac{2\pi i}{2!}\left[\frac{i\xi_n^2+\xi_n-2i}{(\xi+i)^3}\right]^{(2)}\bigg|_{\xi_n=i}dx'\nonumber\\
&=\frac{9h'(0)}{8}g(X_1,X_2)\pi\Omega_3dx'.
\end{align}
Then, when $l=2,$ we have
\begin{align}\label{75}
\Phi_5=\frac{9h'(0)}{8}g(X_1,X_2)\pi\Omega_3dx'.
\end{align}
Similarly, when $l=4,$ we get
\begin{align}\label{19u}
\Phi_5&=-\frac{9h'(0)}{8}[g(X_1,X_2)g(X_3,X_4)-g(X_1,X_3)g(X_2,X_4)+g(X_1,X_4)g(X_2,X_3)]\pi \Omega_3dx'.\nonumber\\
\end{align}
Now $\Phi$ is the sum of the cases (a), (b) and (c), therefore, when $l=2,$ we have
\begin{align}\label{795}
\Phi=\sum_{i=1}^5\Phi_i=-\frac{\partial_{x_n}[g(X_1,X_2)]}{2}\pi\Omega_3dx'.
\end{align}
Similarly, when $l=4,$ we get
\begin{align}\label{109u}
\Phi=\sum_{i=1}^5\Phi_i&=\frac{\partial_{x_n}[g(X_1,X_2)g(X_3,X_4)-g(X_1,X_3)g(X_2,X_4)+g(X_1,X_4)g(X_2,X_3)]}{2}\pi \Omega_3dx'.\nonumber\\
\end{align}
Obviously, when $l=1~or~3,$ we get $\Phi_1=\Phi_2=\Phi_3=\Phi_4=\Phi_5=0.$\\
Then, by (\ref{b15})-(\ref{b1202}), we obtain following theorem
\begin{thm}\label{thmb1}
Let $M$ be a $4$-dimensional oriented
compact manifold with boundary $\partial M$ and the metric
$g^{M}$ be defined as (\ref{b1}), then the following is the generalized noncommutative residue of the Dirac operator\\
(1)when $l=2,$ we get
\begin{align}
\label{b2773}
\widetilde{{\rm Wres}}[\pi^+(LD^{-1})\circ\pi^+D^{-1}]&=32\pi^2\int_{M}\bigg(\frac{1}{3}g(X_1,X_2)s\bigg)d{\rm Vol_{M}}-\int_{\partial M}\bigg(\frac{1}{2}\partial_{x_n}[g(X_1,X_2)]\bigg)\pi\Omega_3d{\rm Vol_{M}},
\end{align}
(2)when $l=4,$ we get
\begin{align}
\label{b12982}
&\widetilde{{\rm Wres}}[\pi^+(LD^{-1})\circ\pi^+D^{-1}]\nonumber\\
&=32\pi^2\int_{M}\bigg(-\frac{1}{3}[g(X_1,X_2)g(X_3,X_4)-g(X_1,X_3)g(X_2,X_4)+g(X_1,X_4)g(X_2,X_3)]s\bigg)d{\rm Vol_{M}}\nonumber\\
&+\int_{\partial M}\bigg(\frac{1}{2}\partial_{x_n}[g(X_1,X_2)g(X_3,X_4)-g(X_1,X_3)g(X_2,X_4)+g(X_1,X_4)g(X_2,X_3)]\bigg)\pi\Omega_3d{\rm Vol_{M}},\nonumber\\
\end{align}
(3)when $l=1~or~3,$ we get
\begin{align}
\label{b1282}
&\widetilde{{\rm Wres}}[\pi^+(LD^{-1})\circ\pi^+D^{-1}]=0.
\end{align}
\end{thm}

\section{A Kastler-Kalau-Walze type theorem for $6$-dimensional manifolds with boundary }
\label{section:4}
Firstly, we prove the Kastler-Kalau-Walze type theorems for the generalized noncommutative
residue $\widetilde{{\rm Wres}}[\pi^+(LD^{-2})\circ\pi^+(D^{-2})]$ on $6$-dimensional manifolds with boundary. From \cite{Wa5}, we know that
\begin{equation}\label{c1}
\widetilde{{\rm Wres}}[\pi^+(LD^{-2})\circ\pi^+(D^{-2})]=\int_M\int_{|\xi'|=1}{\rm
trace}_{\wedge ^*T^*M\bigotimes\mathbb{C}}[\sigma_{-6}(LD^{-4})]\sigma(\xi)dx+\int_{\partial M}\Psi,
\end{equation}
where
\begin{align}\label{c2}
\Psi &=\int_{|\xi'|=1}\int^{+\infty}_{-\infty}\sum^{\infty}_{j, k=0}\sum\frac{(-i)^{|\alpha|+j+k+1}}{\alpha!(j+k+1)!}
\times {\rm trace}_{\wedge ^*T^*M\bigotimes\mathbb{C}}[\partial^j_{x_n}\partial^\alpha_{\xi'}\partial^k_{\xi_n}\sigma^+_{r}(LD^{-2})(x',0,\xi',\xi_n)
\nonumber\\
&\times\partial^\alpha_{x'}\partial^{j+1}_{\xi_n}\partial^k_{x_n}\sigma_{l}
(D^{-2})(x',0,\xi',\xi_n)]d\xi_n\sigma(\xi')dx',
\end{align}
and the sum is taken over $r+\ell-k-j-|\alpha|-1=-6,~~\ r\leq-1,~~\ell\leq -3$.\\
\indent\indent By Theorem \ref{thm2}, we can compute the interior of $\widetilde{{\rm Wres}}[\pi^+(LD^{-2})\circ\pi^+(D^{-2})]$,\\
(1)when $l=2,$ we get
\begin{align}
\label{b146}
&\int_M\int_{|\xi|=1}{\rm
trace}_{\wedge^*T^*M}[\sigma_{-6}(LD^{-4})]\sigma(\xi)dx=128\pi^2\int_{M}\bigg(\frac{2}{3}g(X_1,X_2)s\bigg)d{\rm Vol_{M}},
\end{align}
(2)when $l=4,$ we get
\begin{align}
\label{b12242}
&\int_M\int_{|\xi|=1}{\rm
trace}_{\wedge^*T^*M}[\sigma_{-6}(LD^{-4})]\sigma(\xi)dx\nonumber\\
&=128\pi^2\int_{M}\bigg(-\frac{2}{3}[g(X_1,X_2)g(X_3,X_4)-g(X_1,X_3)g(X_2,X_4)+g(X_1,X_4)g(X_2,X_3)]s\bigg)d{\rm Vol_{M}},
\end{align}
(3)when $l=6,$ we get
\begin{align}
\label{b12342}
&\int_M\int_{|\xi|=1}{\rm
trace}_{\wedge^*T^*M}[\sigma_{-6}(LD^{-4})]\sigma(\xi)dx\nonumber\\
&=128\pi^2\int_{M}\bigg\{-\frac{2}{3}\bigg(g(X_1,X_2)[g(X_3,X_5)g(X_4,X_6)-g(X_3,X_6)g(X_4,X_5)]-g(X_1,X_2)g(X_3,X_4)g(X_5,X_6)\nonumber\\
&+g(X_1,X_3)[g(X_2,X_4)g(X_5,X_6)-g(X_2,X_6)g(X_4,X_5)-g(X_2,X_5)g(X_3,X_6)]+g(X_1,X_4)[g(X_2,X_5)g(X_3,X_6)\nonumber\\
&-g(X_2,X_6)g(X_3,X_5)-g(X_2,X_3)g(X_5,X_6)]+g(X_1,X_5)[g(X_2,X_6)g(X_3,X_4)-g(X_2,X_4)g(X_3,X_6)\nonumber\\
&-g(X_2,X_3)g(X_4,X_6)]+g(X_1,X_6)[g(X_2,X_4)g(X_3,X_5)-g(X_2,X_5)g(X_3,X_4)-g(X_2,X_3)g(X_4,X_5)]\bigg)s\bigg\}d{\rm Vol_{M}},
\end{align}
(4)when $l=1~or~3~or~5,$ we get
\begin{align}
\label{b1pp2}
&&\int_M\int_{|\xi|=1}{\rm
trace}_{\wedge^*T^*M}[\sigma_{-6}(LD^{-4})]\sigma(\xi)dx=0.\nonumber\\
\end{align}
\indent Let the cotangent vector $\xi=\sum_j\xi_jdx_j$, $\xi^j=g^{ij}\xi_i$, $\delta_i=-\frac{1}{4}\sum_{s,t}\omega_{s,t}
(e_i)c(e_s)c(e_t),$ $g^{ij}=g(dx_i,dx_j)$ and $\nabla^L_{\partial_i}\partial_j=\sum_k\Gamma^k_{ij}\partial_k,$ $\Gamma^k=g^{ij}\Gamma^k_{ij}$, $\delta^j=g^{ij}\delta_i.$ Then, by \cite{Ka}, we obtain
\begin{lem}\label{clem2} The following identities hold:
\begin{align}\label{c10}
\sigma_{-2}(LD^{-2})&=L|\xi|^{-2};\nonumber\\
\sigma_{-3}(LD^{-2})&=-\sqrt{-1}|\xi|^{-4}L\xi_k(\Gamma^k-2\delta^k)-2\sqrt{-1}|\xi|^{-6}L\xi^j\xi_\alpha\xi_\beta\partial_jg^{\alpha\beta}.\nonumber\\
\end{align}
\end{lem}
When $n=6$, then ${\rm tr}_{S(TM)}[\texttt{id}]=8$.
Since the sum is taken over $-r-\ell+k+j+|\alpha|-1=-6, \ r, \ell\leq -2$, then we have the $\int_{\partial_{M}}\Psi$
is the sum of the following five cases:\\
\noindent  {\bf case (a)~(I)}~$r=-2,~\ell=-2,~j=k=0,~|\alpha|=1$.\\
By (\ref{c2}), we get
 \begin{equation}\label{c12}
\Psi_1=-\int_{|\xi'|=1}\int^{+\infty}_{-\infty}\sum_{|\alpha|=1}{\rm trace}
\Big[\partial^{\alpha}_{\xi'}\pi^{+}_{\xi_{n}}\sigma_{-1}(LD^{-2})
      \times\partial^{\alpha}_{x'}\partial_{\xi_{n}}\sigma_{-2}(D^{-2})\Big](x_0)d\xi_n\sigma(\xi')dx'.
\end{equation}
By Lemma \ref{clem2}, for $i<n$, we have
 \begin{align}\label{c13}
 \partial_{x_{i}}\sigma_{-2}(D^{-2})(x_0)&=
      \partial_{x_{i}}(|\xi|^{-2})(x_{0})\nonumber\\
      &=-\frac{\partial_{x_{i}}(|\xi|^{2})(x_{0})}{|\xi|^{4}}\nonumber\\
      &=0,\nonumber\\
\end{align}
\noindent so when $l$ is an integer between 1 and 6, we have $\Psi_1=0$.\\
~\\
\noindent  {\bf case (a)~(II)}~$r=-2,~\ell=-2,~|\alpha|=k=0,~j=1$.\\
By (\ref{c2}), we have
  \begin{equation}\label{c14}
\Psi_2=-\frac{1}{2}\int_{|\xi'|=1}\int^{+\infty}_{-\infty} {\rm
trace} \Big[\partial_{x_{n}}\pi^{+}_{\xi_{n}}\sigma_{-2}(LD^{-2})
      \times\partial^{2}_{\xi_{n}}\sigma_{-2}(D^{-2})\Big](x_0)d\xi_n\sigma(\xi')dx'.
\end{equation}
By computations, we have\\
\begin{align}\label{c15}
\partial^{2}_{\xi_{n}}\sigma_{-2}(D^{-2})(x_0)&=\partial^{2}_{\xi_{n}}(|\xi|^{-2})(x_0)=2\frac{3\xi_n^2-1}{(1+\xi_n^2)^3}.\nonumber\\
\end{align}
\begin{align}\label{c115}
\partial_{x_{n}}\sigma_{-2}(LD^{-2})(x_0)|_{|\xi'|=1}&=\frac{\partial_{x_{n}}[\sum_{\alpha,\beta=1}^6a_{1\alpha}a_{2\beta}c(e_\alpha)c(e_\beta)]}{|\xi|^2}-\frac{\sum_{\alpha,\beta=1}^6a_{1\alpha}a_{2\beta}c(e_\alpha)c(e_\beta)h'(0)}{|\xi|^4}\nonumber\\
&=\frac{\sum_{\alpha,\beta=1}^6\partial_{x_{n}}(a_{1\alpha})a_{2\beta}c(e_\alpha)c(e_\beta)}{|\xi|^2}+\frac{\sum_{\alpha,\beta=1}^6a_{1\alpha}\partial_{x_{n}}(a_{2\beta})c(e_\alpha)c(e_\beta)}{|\xi|^2}\nonumber\\
&-\frac{\sum_{\alpha,\beta=1}^6a_{1\alpha}a_{2\beta}c(e_\alpha)c(e_\beta)h'(0)}{|\xi|^4}.\nonumber\\
\end{align}
\begin{align}\label{c1995}
&\pi^+_{\xi_n}\partial_{x_{n}}\sigma_{-2}(LD^{-2})(x_0)|_{|\xi'|=1}\nonumber\\
&=-\frac{i}{2(\xi_n-i)}\sum_{\alpha,\beta=1}^6\partial_{x_{n}}(a_{1\alpha})a_{2\beta}c(e_\alpha)c(e_\beta)-\frac{i}{2(\xi_n-i)}\sum_{\alpha,\beta=1}^6a_{1\alpha}\partial_{x_{n}}(a_{2\beta})c(e_\alpha)c(e_\beta)\nonumber\\
&+\frac{(i\xi_n+2)h'(0)}{4(\xi_n-i)^2}\sum_{\alpha,\beta=1}^6a_{1\alpha}a_{2\beta}c(e_\alpha)c(e_\beta).\nonumber\\
\end{align}
Since $n=6$, ${\rm tr}[-\texttt{id}]=-8$. By the relation of the Clifford action and ${\rm tr}ab={\rm tr}ba$,  then
\begin{align}\label{c16}
&{\rm tr}[c(\xi')c(dx_n)]=0;~~{\rm tr}[c(dx_n)^2]=-8;~~{\rm tr}[c(\xi')^2](x_0)|_{|\xi'|=1}=-8;~~{\rm tr}[\partial_{x_n}c(\xi')c(dx_n)]=0;\nonumber\\
&{\rm tr}[\partial_{x_n}c(\xi')c(\xi')](x_0)|_{|\xi'|=1}=-4h'(0);~~\sum_{\alpha,\beta=1}^6\partial_{x_n}(a_{1\alpha})(a_{2\alpha})+\sum_{\alpha,\beta=1}^6a_{1\alpha}\partial_{x_n}(a_{2\alpha})=\partial_{x_n}[g(X_1,X_2)];\nonumber\\
&{\rm tr}[\sum_{\alpha,\beta=1}^6a_{1\alpha}a_{2\beta}c(e_\alpha)c(e_\beta)c(\xi')c(\xi')]={\rm tr}[\sum_{\alpha,\beta=1}^6a_{1\alpha}a_{2\beta}c(e_\alpha)c(e_\beta)c(dx_n)c(dx_n)]=8g(X_1,X_2);\nonumber\\
&{\rm tr}[\sum_{\alpha,\beta=1}^6\partial_{x_n}(a_{1\alpha})a_{2\beta}c(e_\alpha)c(e_\beta)c(\xi')c(\xi')]={\rm tr}[\sum_{\alpha,\beta=1}^6\partial_{x_n}(a_{1\alpha})a_{2\beta}c(e_\alpha)c(e_\beta)c(dx_n)c(dx_n)]=8\sum_{\alpha,\beta=1}^6\partial_{x_n}(a_{1\alpha})a_{2\alpha}.\nonumber\\
\end{align}
When $l=2,$ we have
\begin{align}\label{871}
&{\rm
trace} \Big[\partial_{x_{n}}\pi^{+}_{\xi_{n}}\sigma_{-2}(LD^{-2})
      \times\partial^{2}_{\xi_{n}}\sigma_{-2}(D^{-2})\Big](x_0)\nonumber\\
&=8\frac{i(3\xi_n^2-1)}{(\xi_n-i)^4(\xi_n+i)^3}\partial_{x_{n}}[g(X_1,X_2)]-4h'(0)\frac{(i\xi_n+2)(3\xi_n^2-1)}{(\xi_n-i)^5(\xi_n+i)^3}g(X_1,X_2).
\end{align}
Then, we obtain
\begin{align}\label{c18}
\Psi_2&=-\frac{1}{2}\int_{|\xi'|=1}\int^{+\infty}_{-\infty} \bigg\{8\frac{i(3\xi_n^2-1)}{(\xi_n-i)^4(\xi_n+i)^3}\partial_{x_{n}}[g(X_1,X_2)]-4h'(0)\frac{(i\xi_n+2)(3\xi_n^2-1)}{(\xi_n-i)^5(\xi_n+i)^3}g(X_1,X_2)\bigg\}d\xi_n\sigma(\xi')dx'\nonumber\\
     &=-4i\partial_{x_{n}}[g(X_1,X_2)]\Omega_{4}\int_{\Gamma^{+}}\frac{3\xi_n^2-1}{(\xi_n-i)^4(\xi_n+i)^3}d\xi_{n}dx'+2h'(0)g(X_1,X_2)\Omega_{4}\int_{\Gamma^{+}}\frac{(i\xi_n+2)(3\xi_n^2-1)}{(\xi_n-i)^5(\xi_n+i)^3}d\xi_{n}dx'\nonumber\\
     &=-4i\partial_{x_{n}}[g(X_1,X_2)]\Omega_{4}\frac{2\pi i}{3!}\left[\frac{3\xi_n^2-1}{(\xi_n+i)^3}\right]^{(3)}\bigg|_{|\xi'|=1}dx'+2h'(0)g(X_1,X_2)\Omega_{4}\frac{2\pi i}{4!}\left[\frac{(i\xi_n+2)(3\xi_n^2-1)}{(\xi_n+i)^3}\right]^{(4)}\bigg|_{|\xi'|=1}dx'\nonumber\\
     &=\left(\frac{5h'(0)}{8}g(X_1,X_2)-\partial_{x_{n}}[g(X_1,X_2)]\right)\pi\Omega_{4}dx',
\end{align}
where ${\rm \Omega_{4}}$ is the canonical volume of $S^{4}.$\\
Similarly, when $l=4,$ we get
\begin{align}\label{1um}
\Psi_2&=\bigg(-\frac{5h'(0)}{8}[g(X_1,X_2)g(X_3,X_4)-g(X_1,X_3)g(X_2,X_4)+g(X_1,X_4)g(X_2,X_3)]\nonumber\\
&+\partial_{x_n}[g(X_1,X_2)g(X_3,X_4)-g(X_1,X_3)g(X_2,X_4)+g(X_1,X_4)g(X_2,X_3)]\bigg)\pi\Omega_{4}dx'.\nonumber\\
\end{align}
Similarly, when $l=6,$ we get
\begin{align}\label{1bbm}
\Psi_2&=\bigg\{-\frac{5h'(0)}{8}\bigg(g(X_1,X_2)[g(X_3,X_5)g(X_4,X_6)-g(X_3,X_6)g(X_4,X_5)-g(X_3,X_4)g(X_5,X_6)]\nonumber\\
&+g(X_1,X_3)[g(X_2,X_4)g(X_5,X_6)-g(X_2,X_6)g(X_4,X_5)-g(X_2,X_5)g(X_3,X_6)]+g(X_1,X_4)[g(X_2,X_5)g(X_3,X_6)\nonumber\\
&-g(X_2,X_6)g(X_3,X_5)-g(X_2,X_3)g(X_5,X_6)]+g(X_1,X_5)[g(X_2,X_6)g(X_3,X_4)-g(X_2,X_4)g(X_3,X_6)\nonumber\\
&-g(X_2,X_3)g(X_4,X_6)]+g(X_1,X_6)[g(X_2,X_4)g(X_3,X_5)-g(X_2,X_5)g(X_3,X_4)-g(X_2,X_3)g(X_4,X_5)]\bigg)\nonumber\\
&+\partial_{x_n}\bigg(g(X_1,X_2)[g(X_3,X_5)g(X_4,X_6)-g(X_3,X_6)g(X_4,X_5)-g(X_3,X_4)g(X_5,X_6)]\nonumber\\
&+g(X_1,X_3)[g(X_2,X_4)g(X_5,X_6)-g(X_2,X_6)g(X_4,X_5)-g(X_2,X_5)g(X_3,X_6)]+g(X_1,X_4)[g(X_2,X_5)g(X_3,X_6)\nonumber\\
&-g(X_2,X_6)g(X_3,X_5)-g(X_2,X_3)g(X_5,X_6)]+g(X_1,X_5)[g(X_2,X_6)g(X_3,X_4)-g(X_2,X_4)g(X_3,X_6)\nonumber\\
&-g(X_2,X_3)g(X_4,X_6)]+g(X_1,X_6)[g(X_2,X_4)g(X_3,X_5)-g(X_2,X_5)g(X_3,X_4)-g(X_2,X_3)g(X_4,X_5)]\bigg)\bigg\}\pi\Omega_{4}dx'.\nonumber\\
\end{align}
\noindent  {\bf case (a)~(III)}~$r=-2,~\ell=-2,~|\alpha|=j=0,~k=1$.\\
By (\ref{c2}), we have
 \begin{align}\label{c19}
\Psi_3&=-\frac{1}{2}\int_{|\xi'|=1}\int^{+\infty}_{-\infty}{\rm trace} \Big[\partial_{\xi_{n}}\pi^{+}_{\xi_{n}}\sigma_{-2}(LD^{-2})
      \times\partial_{\xi_{n}}\partial_{x_{n}}\sigma_{-2}(D^{-2})\Big](x_0)d\xi_n\sigma(\xi')dx'\nonumber\\
      &=\frac{1}{2}\int_{|\xi'|=1}\int^{+\infty}_{-\infty}{\rm trace} \Big[\partial_{\xi_{n}}^2\pi^{+}_{\xi_{n}}\sigma_{-2}(LD^{-2})
      \times\partial_{x_{n}}\sigma_{-2}(D^{-2})\Big](x_0)d\xi_n\sigma(\xi')dx'\nonumber\\
\end{align}
By computations, when $l=2,$ we have\\
\begin{equation}\label{c20}
\partial_{x_{n}}\sigma_{-2}(D^{-2})(x_0)|_{|\xi'|=1}=-\frac{h'(0)}{(1+\xi_{n}^{2})^2}.
\end{equation}
\begin{align}\label{c520}
\partial_{\xi_{n}}^2\pi^{+}_{\xi_{n}}\sigma_{-2}(LD^{-2})(x_0)|_{|\xi'|=1}=-\frac{\sum_{\alpha,\beta=1}^6a_{1\alpha}a_{2\beta}c(e_\alpha)c(e_\beta)i}{(\xi_{n}-i)^3}.
\end{align}
Combining (\ref{c1995}) and (\ref{c20}), we have
\begin{equation}\label{c21}
{\rm trace} \Big[\partial_{\xi_{n}}^2\pi^{+}_{\xi_{n}}\sigma_{-2}(LD^{-2})
      \times\partial_{x_{n}}\sigma_{-2}(D^{-2})\Big](x_0)(x_{0})|_{|\xi'|=1}
=-\frac{8h'(0)i}{(\xi_{n}-i)^{5}(\xi+i)^2}g(X_1,X_2).
\end{equation}
Then
\begin{align}\label{c22}
\Psi_3&=\frac{1}{2}\int_{|\xi'|=1}\int^{+\infty}_{-\infty} -\frac{8h'(0)i}{(\xi_{n}-i)^{5}(\xi+i)^2}g(X_1,X_2)d\xi_n\sigma(\xi')dx'\nonumber\\
     &=-4h'(0)ig(X_1,X_2)\Omega_{4}\int_{\Gamma^{+}}\frac{1}{(\xi_{n}-i)^{5}(\xi+i)^2}d\xi_{n}dx'\nonumber\\
     &=-4h'(0)ig(X_1,X_2)\Omega_{4}\frac{2\pi i}{4!}\left[\frac{1}{(\xi+i)^2}\right]^{(4)}\bigg|_{|\xi'|=1}dx'\nonumber\\
     &=\frac{-5h'(0)}{8}g(X_1,X_2)\pi\Omega_{4}dx'.
\end{align}
Similarly, when $l=4,$ we get
\begin{align}\label{1eum}
\Psi_3&=\frac{5h'(0)}{8}[g(X_1,X_2)g(X_3,X_4)-g(X_1,X_3)g(X_2,X_4)+g(X_1,X_4)g(X_2,X_3)]\pi\Omega_{4}dx'.\nonumber\\
\end{align}
When $l=6,$ we get
\begin{align}\label{1elum}
\Psi_3&=\frac{5h'(0)}{8}\bigg\{g(X_1,X_2)[g(X_3,X_5)g(X_4,X_6)-g(X_3,X_6)g(X_4,X_5)-g(X_3,X_4)g(X_5,X_6)]\nonumber\\
&+g(X_1,X_3)[g(X_2,X_4)g(X_5,X_6)-g(X_2,X_6)g(X_4,X_5)-g(X_2,X_5)g(X_3,X_6)]\nonumber\\
&+g(X_1,X_4)[g(X_2,X_5)g(X_3,X_6)-g(X_2,X_6)g(X_3,X_5)-g(X_2,X_3)g(X_5,X_6)]\nonumber\\
&+g(X_1,X_5)[g(X_2,X_6)g(X_3,X_4)-g(X_2,X_4)g(X_3,X_6)-g(X_2,X_3)g(X_4,X_6)]\nonumber\\
&+g(X_1,X_6)[g(X_2,X_4)g(X_3,X_5)-g(X_2,X_5)g(X_3,X_4)-g(X_2,X_3)g(X_4,X_5)]\bigg\}\pi\Omega_{4}dx'.\nonumber\\
\end{align}
\noindent  {\bf case (b)}~$r=-2,~\ell=-3,~|\alpha|=j=k=0$.\\
By (\ref{c2}), we have
 \begin{align}\label{c23}
\Psi_4&=-i\int_{|\xi'|=1}\int^{+\infty}_{-\infty}{\rm trace} \Big[\pi^{+}_{\xi_{n}}\sigma_{-2}(LD^{-2})
      \times\partial_{\xi_{n}}\sigma_{-3}(D^{-2})\Big](x_0)d\xi_n\sigma(\xi')dx'\nonumber\\
&=i\int_{|\xi'|=1}\int^{+\infty}_{-\infty}{\rm trace} [\partial_{\xi_n}\pi^+_{\xi_n}\sigma_{-2}(LD^{-2})\times
\sigma_{-3}(D^{-2})](x_0)d\xi_n\sigma(\xi')dx'.
\end{align}
When $l=2,$ we have
\begin{align}\label{c000}
\partial_{\xi_{n}}\pi^+_{\xi_n}\sigma_{-2}(LD^{-2})(x_0)|_{|\xi'|=1}=\frac{\sum_{\alpha,\beta=1}^6a_{1\alpha}c(e_\alpha)a_{2\beta}c(e_\beta)i}{2(\xi_{n}-i)^2}.
\end{align}

In the normal coordinate, $g^{ij}(x_{0})=\delta^{j}_{i}$ and $\partial_{x_{j}}(g^{\alpha\beta})(x_{0})=0$, if $j<n$; $\partial_{x_{j}}(g^{\alpha\beta})(x_{0})=h'(0)\delta^{\alpha}_{\beta}$, if $j=n$.
So by  \cite{Wa3}, when $k<n$, we have $\Gamma^{n}(x_{0})=\frac{5}{2}h'(0)$, $\Gamma^{k}(x_{0})=0$, $\delta^{n}(x_{0})=0$ and $\delta^{k}=\frac{1}{4}h'(0)c(e_{k})c(e_{n})$. Then, we obtain

\begin{align}\label{c24}
\sigma_{-3}(D^{-2})(x_{0})|_{|\xi'|=1}&=
-\sqrt{-1}|\xi|^{-4}\xi_k(\Gamma^k-2\delta^k)(x_{0})|_{|\xi'|=1}-\sqrt{-1}|\xi|^{-6}2\xi^j\xi_\alpha\xi_\beta\partial_jg^{\alpha\beta}(x_{0})|_{|\xi'|=1}\nonumber\\
&=\frac{-i}{(1+\xi_n^2)^2}\Big(-\frac{1}{2}h'(0)\sum_{k<n}\xi_k
\widetilde{c}(e_k)\widetilde{c}(e_n)+\frac{5}{2}h'(0)\xi_n\Big)-\frac{2ih'(0)\xi_n}{(1+\xi_n^2)^3}.\nonumber\\
\end{align}
We note that $i<n,~\int_{|\xi'|=1}\{\xi_{i_{1}}\xi_{i_{2}}\cdots\xi_{i_{2d+1}}\}\sigma(\xi')=0$, so the first term in (\ref{c24})  has no contribution for computing {\bf case (b)}.\\
By (\ref{c000}) and (\ref{c24}), we have
\begin{align}\label{c25}
&{\rm trace} [\partial_{\xi_n}\pi^+_{\xi_n}\sigma_{-2}(LD^{-2})\times
\sigma_{-3}(D^{-2})](x_0)|_{|\xi'|=1} \nonumber\\
&=-\frac{2h'(0)\xi_n(5\xi_n^2-1)}{(\xi_n-i)^5(\xi_n+i)^3}g(X_1,X_2).
\end{align}
So when $l=2$, we have
\begin{align}\label{c28}
\Psi_4&=
 i\int_{|\xi'|=1}\int^{+\infty}_{-\infty}-\frac{2h'(0)\xi_n(5\xi_n^2-1)}{(\xi_n-i)^5(\xi_n+i)^3}g(X_1,X_2)d\xi_n\sigma(\xi')dx'\nonumber\\
&=-2ih'(0)g(X_1,X_2)\Omega_{4}\int_{\Gamma^{+}}\frac{\xi_n(5\xi_n^2-1)}{(\xi_n-i)^5(\xi_n+i)^3}d\xi_{n}dx'\nonumber\\
&=-2ih'(0)g(X_1,X_2)\Omega_{4}\frac{2\pi i}{4!}\left[\frac{\xi_n(5\xi_n^2-1)}{(\xi_n+i)^3}\right]^{(4)}\bigg|_{|\xi'|=1}dx'\nonumber\\
&=\frac{15h'(0)}{8}g(X_1,X_2)\pi\Omega_4dx'.\nonumber\\
\end{align}
Similarly, when $l=4,$ we get
\begin{align}\label{1eubbm}
\Psi_4&=-\frac{15h'(0)}{8}[g(X_1,X_2)g(X_3,X_4)-g(X_1,X_3)g(X_2,X_4)+g(X_1,X_4)g(X_2,X_3)]\pi\Omega_{4}dx'.\nonumber\\
\end{align}
When $l=6,$ we get
\begin{align}\label{1eglum}
\Psi_4&=-\frac{15h'(0)}{8}\bigg\{g(X_1,X_2)[g(X_3,X_5)g(X_4,X_6)-g(X_3,X_6)g(X_4,X_5)-g(X_3,X_4)g(X_5,X_6)]\nonumber\\
&+g(X_1,X_3)[g(X_2,X_4)g(X_5,X_6)-g(X_2,X_6)g(X_4,X_5)-g(X_2,X_5)g(X_3,X_6)]\nonumber\\
&+g(X_1,X_4)[g(X_2,X_5)g(X_3,X_6)-g(X_2,X_6)g(X_3,X_5)-g(X_2,X_3)g(X_5,X_6)]\nonumber\\
&+g(X_1,X_5)[g(X_2,X_6)g(X_3,X_4)-g(X_2,X_4)g(X_3,X_6)-g(X_2,X_3)g(X_4,X_6)]\nonumber\\
&+g(X_1,X_6)[g(X_2,X_4)g(X_3,X_5)-g(X_2,X_5)g(X_3,X_4)-g(X_2,X_3)g(X_4,X_5)]\bigg\}\pi\Omega_{4}dx'.\nonumber\\
\end{align}
\noindent {\bf  case (c)}~$r=-3,~\ell=-2,~|\alpha|=j=k=0$.\\
By (\ref{c2}), we have
\begin{align}\label{c29}
\Psi_5&=-i\int_{|\xi'|=1}\int^{+\infty}_{-\infty}{\rm trace} \Big[\pi^{+}_{\xi_{n}}\sigma_{-3}(LD^{-2})
      \times\partial_{\xi_{n}}\sigma_{-2}(D^{-2})\Big](x_0)d\xi_n\sigma(\xi')dx'\nonumber\\
      &={\bf  case (b)}-i\int_{|\xi'|=1}\int^{+\infty}_{-\infty}{\rm trace} \Big[\pi^{+}_{\xi_{n}}\partial_{\xi_{n}}\sigma_{-2}(D^{-2})
      \times\sigma_{-3}(LD^{-2})\Big](x_0)d\xi_n\sigma(\xi')dx'.\nonumber\\
\end{align}
Then, we only need to compute\\
\begin{align}\label{c429}
-i\int_{|\xi'|=1}\int^{+\infty}_{-\infty}{\rm trace} \Big[\pi^{+}_{\xi_{n}}\partial_{\xi_{n}}\sigma_{-2}(D^{-2})
      \times\sigma_{-3}(LD^{-2})\Big](x_0)d\xi_n\sigma(\xi')dx'.\nonumber\\
\end{align}
By Lemma \ref{clem2}, when $l=2,$ we have
\begin{align}\label{c30}
\partial_{\xi_n}\sigma_{-2}(D^{-2})(x_0)&=-\frac{2\xi_n}{(1+\xi_n^2)^2}.\nonumber\\
\end{align}
\begin{align}\label{c31}
\sigma_{-3}(LD^{-2})=\sum^6_{\alpha,\beta=1}a_{1\alpha}a_{2\beta}c(e_\alpha)c(e_\beta)\bigg(\frac{-i}{(1+\xi_n^2)^2}\Big(-\frac{1}{2}h'(0)\sum_{k<n}\xi_k
\widetilde{c}(e_k)\widetilde{c}(e_n)+\frac{5}{2}h'(0)\xi_n\Big)-\frac{2ih'(0)\xi_n}{(1+\xi_n^2)^3}\bigg).\nonumber\\
\end{align}
We note that $i<n,~\int_{|\xi'|=1}\{\xi_{i_{1}}\xi_{i_{2}}\cdots\xi_{i_{2d+1}}\}\sigma(\xi')=0$, so the first term in (\ref{c31})  has no contribution for computing {\bf case (c)}.\\
On the other hand,
\begin{align}\label{c32}
{\rm trace} \Big[\pi^{+}_{\xi_{n}}\partial_{\xi_{n}}\sigma_{-2}(D^{-2})
      \times\sigma_{-3}(LD^{-2})\Big](x_0)&=-8\frac{i h'(0)\xi_{n}^2(9+5\xi_n^2)}{(1+\xi_{n}^{2})^5}g(X_1,X_2).\nonumber\\
\end{align}
Then, when $l=2,$ we have
\begin{align}\label{c45}
\Psi_5&=\frac{15h'(0)}{8}g(X_1,X_2)\pi\Omega_4dx'
 -i\int_{|\xi'|=1}\int^{+\infty}_{-\infty}
 -8\frac{i h'(0)\xi_{n}^2(9+5\xi_n^2)}{(1+\xi_{n}^{2})^5}g(X_1,X_2)d\xi_n\sigma(\xi')dx' \nonumber\\
&=\frac{15h'(0)}{8}g(X_1,X_2)\pi\Omega_4dx'-8h'(0)g(X_1,X_2)\Omega_4\int_{\Gamma^{+}}\frac{\xi_{n}^2(9+5\xi_n^2)}{(1+\xi_{n}^{2})^5}d\xi_{n}dx'\nonumber\\
&=\frac{15h'(0)}{8}g(X_1,X_2)\pi\Omega_4dx'-8h'(0)g(X_1,X_2)\Omega_4\frac{2\pi i}{4!}\left[\frac{\xi_{n}^2(9+5\xi_n^2)}{(\xi_n+i)^5}\right]^{(4)}\bigg|_{|\xi'|=1}dx'\nonumber\\
&=-\frac{15h'(0)}{8}g(X_1,X_2)\pi\Omega_4dx'.\nonumber\\
\end{align}
Similarly, when $l=4,$ we get
\begin{align}\label{12eum}
\Psi_5&=\frac{15h'(0)}{8}[g(X_1,X_2)g(X_3,X_4)-g(X_1,X_3)g(X_2,X_4)+g(X_1,X_4)g(X_2,X_3)]\pi\Omega_{4}dx'.\nonumber\\
\end{align}
When $l=6,$ we get
\begin{align}\label{12elum}
\Psi_5&=\frac{15h'(0)}{8}\bigg\{g(X_1,X_2)[g(X_3,X_5)g(X_4,X_6)-g(X_3,X_6)g(X_4,X_5)-g(X_3,X_4)g(X_5,X_6)]\nonumber\\
&+g(X_1,X_3)[g(X_2,X_4)g(X_5,X_6)-g(X_2,X_6)g(X_4,X_5)-g(X_2,X_5)g(X_3,X_6)]\nonumber\\
&+g(X_1,X_4)[g(X_2,X_5)g(X_3,X_6)-g(X_2,X_6)g(X_3,X_5)-g(X_2,X_3)g(X_5,X_6)]\nonumber\\
&+g(X_1,X_5)[g(X_2,X_6)g(X_3,X_4)-g(X_2,X_4)g(X_3,X_6)-g(X_2,X_3)g(X_4,X_6)]\nonumber\\
&+g(X_1,X_6)[g(X_2,X_4)g(X_3,X_5)-g(X_2,X_5)g(X_3,X_4)-g(X_2,X_3)g(X_4,X_5)]\bigg\}\pi\Omega_{4}dx'.\nonumber\\
\end{align}
Now $\Psi$ is the sum of the cases (a), (b) and (c), then when $l=2,$ we get
\begin{align}\label{795}
\Psi=\sum_{i=1}^5\Psi_i=-\partial_{x_n}[g(X_1,X_2)]\pi\Omega_4dx'.
\end{align}
Similarly, when $l=4,$ we get
\begin{align}\label{109u}
\Psi=\sum_{i=1}^5\Psi_i&=\partial_{x_n}[g(X_1,X_2)g(X_3,X_4)-g(X_1,X_3)g(X_2,X_4)+g(X_1,X_4)g(X_2,X_3)]\pi \Omega_4dx'.\nonumber\\
\end{align}
Similarly, when $l=6,$ we get
\begin{align}\label{119u}
\Psi=\sum_{i=1}^5\Psi_i&=\partial_{x_n}\bigg(g(X_1,X_2)[g(X_3,X_5)g(X_4,X_6)-g(X_3,X_6)g(X_4,X_5)-g(X_3,X_4)g(X_5,X_6)]\nonumber\\
&+g(X_1,X_3)[g(X_2,X_4)g(X_5,X_6)-g(X_2,X_6)g(X_4,X_5)-g(X_2,X_5)g(X_3,X_6)]+g(X_1,X_4)\nonumber\\
&[g(X_2,X_5)(X_3,X_6)-g(X_2,X_6)g(X_3,X_5)-g(X_2,X_3)g(X_5,X_6)]+g(X_1,X_5)[g(X_2,X_6)\nonumber\\
&g(X_3,X_4)-g(X_2,X_4)g(X_3,X_6)-g(X_2,X_3)g(X_4,X_6)]+g(X_1,X_6)[g(X_2,X_4)g(X_3,X_5)\nonumber\\
&-g(X_2,X_5)g(X_3,X_4)-g(X_2,X_3)g(X_4,X_5)]\bigg)\pi \Omega_4dx'.\nonumber\\
\end{align}
Obviously, when $l=1~or~3~or~5,$ we get $\Psi_1=\Psi_2=\Psi_3=\Psi_4=\Psi_5=0.$\\
By (\ref{c1})-(\ref{b1pp2}), we obtain following theorem
\begin{thm}\label{cthm1}
Let $M$ be a $6$-dimensional oriented
compact manifold with boundary $\partial M$ and the metric
$g^{M}$ be defined as (\ref{b1}), then the following are the generalized noncommutative residue of the Dirac operator\\
(1)when $l=2,$ we get
\begin{align}
\label{b263}
\widetilde{{\rm Wres}}[\pi^+(LD^{-2})\circ\pi^+(D^{-2})]&=128\pi^2\int_{M}\bigg(\frac{2}{3}g(X_1,X_2)s\bigg)d{\rm Vol_{M}}+\int_{\partial M}\bigg(-\partial_{x_n}[g(X_1,X_2)]\bigg)\pi\Omega_4d{\rm Vol_{M}}.
\end{align}
(2)when $l=4,$ we get
\begin{align}
\label{b12982}
&\widetilde{{\rm Wres}}[\pi^+(LD^{-2})\circ\pi^+(D^{-2})]\nonumber\\
&=128\pi^2\int_{M}\bigg(-\frac{2}{3}[g(X_1,X_2)g(X_3,X_4)-g(X_1,X_3)g(X_2,X_4)+g(X_1,X_4)g(X_2,X_3)]s\bigg)d{\rm Vol_{M}}\nonumber\\
&+\int_{\partial M}\bigg(\partial_{x_n}[g(X_1,X_2)g(X_3,X_4)-g(X_1,X_3)g(X_2,X_4)+g(X_1,X_4)g(X_2,X_3)]\bigg)\pi \Omega_4d{\rm Vol_{M}}.\nonumber\\
\end{align}
(3)when $l=6,$ we get
\begin{align}
\label{b1q42}
&\widetilde{{\rm Wres}}[\pi^+(LD^{-2})\circ\pi^+(D^{-2})]\nonumber\\
&=128\pi^2\int_{M}\bigg\{-\frac{2}{3}\bigg(g(X_1,X_2)[g(X_3,X_5)g(X_4,X_6)-g(X_3,X_6)g(X_4,X_5)-g(X_3,X_4)g(X_5,X_6)]\nonumber\\
&+g(X_1,X_3)[g(X_2,X_4)g(X_5,X_6)-g(X_2,X_6)g(X_4,X_5)-g(X_2,X_5)g(X_3,X_6)]+g(X_1,X_4)[g(X_2,X_5)g(X_3,X_6)\nonumber\\
&-g(X_2,X_6)g(X_3,X_5)-g(X_2,X_3)g(X_5,X_6)]+g(X_1,X_5)[g(X_2,X_6)g(X_3,X_4)-g(X_2,X_4)g(X_3,X_6)\nonumber\\
&-g(X_2,X_3)g(X_4,X_6)]+g(X_1,X_6)[g(X_2,X_4)g(X_3,X_5)-g(X_2,X_5)g(X_3,X_4)-g(X_2,X_3)g(X_4,X_5)]\bigg)s\bigg\}d{\rm Vol_{M}}\nonumber\\
&+\int_{\partial M}\bigg\{\partial_{x_n}\bigg(g(X_1,X_2)[g(X_3,X_5)g(X_4,X_6)-g(X_3,X_6)g(X_4,X_5)-g(X_3,X_4)g(X_5,X_6)]\nonumber\\
&+g(X_1,X_3)[g(X_2,X_4)g(X_5,X_6)-g(X_2,X_6)g(X_4,X_5)-g(X_2,X_5)g(X_3,X_6)]+g(X_1,X_4)[g(X_2,X_5)g(X_3,X_6)\nonumber\\
&-g(X_2,X_6)g(X_3,X_5)-g(X_2,X_3)g(X_5,X_6)]+g(X_1,X_5)[g(X_2,X_6)g(X_3,X_4)-g(X_2,X_4)g(X_3,X_6)\nonumber\\
&-g(X_2,X_3)g(X_4,X_6)]+g(X_1,X_6)[g(X_2,X_4)g(X_3,X_5)-g(X_2,X_5)g(X_3,X_4)-g(X_2,X_3)g(X_4,X_5)]\bigg)\bigg\}\pi\Omega_4d{\rm Vol_{M}}.\nonumber\\
\end{align}
(4)when $l=1~or~3~or~5,$ we get
\begin{align}
\label{b1442}
&\widetilde{{\rm Wres}}[\pi^+(LD^{-2})\circ\pi^+(D^{-2})]=0.
\end{align}
\end{thm}

\section*{Acknowledgements}
This work was supported by NSFC. 11771070 .
 The authors thank the referee for his (or her) careful reading and helpful comments.

\end{document}